\documentclass[english]{bourbaki}
\usepackage{times}

\usepackage{amsthm}
\usepackage{amscd}
\usepackage{amsfonts}
\usepackage{amssymb}

\theoremstyle{plain}
\newtheorem{cor}{Corollary}
\newtheorem{thm}{Theorem}
\newtheorem{propn}{Proposition}
\newtheorem{lem}{Lemma}
\newtheorem{defn}{Definition}

\theoremstyle{definition}

\theoremstyle{remark}
\newtheorem{rem}{Remark}

\newcommand{\beq}{\begin{equation}}
\newcommand{\eeq}{\end{equation}}

\newcommand{\id}{\mbox{id}}

\newcommand{\pa}{\partial}
\newcommand{\ot}{\otimes}
\newcommand{\ra}{\rightarrow}

\newcommand{\ti}{\times}

\newcommand{\fr}[2]{{\textstyle \frac{#1}{#2} }}

\newcommand{\fsl}{{\mathfrak s}{\mathfrak l}}

\newcommand{\bra}{\langle}
\newcommand{\ket}{\rangle}
\newcommand{\al}{\alpha}
\newcommand{\bal}{\bar{\alpha}}
\newcommand{\be}{\beta}
\newcommand{\ga}{\gamma}

\newcommand{\Ga}{\Gamma}
\newcommand{\de}{\delta}
\newcommand{\De}{\Delta}
\newcommand{\ep}{\epsilon}
\newcommand{\ka}{\kappa}
\newcommand{\la}{\lambda}
\newcommand{\om}{\omega}

\newcommand{\si}{\sigma}
\newcommand{\up}{\Upsilon}

\newcommand{\TB}{\tilde{B}}
\newcommand{\tE}{\tilde{E}}
\newcommand{\tF}{\tilde{F}}
\newcommand{\tK}{\tilde{K}}

\newcommand{\CA}{{\mathcal A}}
\newcommand{\CB}{{\mathcal B}}
\newcommand{\CC}{{\mathcal C}}
\newcommand{\cD}{{\mathcal D}}

\newcommand{\CF}{{\mathcal F}}

\newcommand{\CH}{{\mathcal H}}
\newcommand{\CI}{{\mathcal I}}

\newcommand{\CK}{{\mathcal K}}
\newcommand{\CL}{{\mathcal L}}

\newcommand{\CO}{{\mathcal O}}  
\newcommand{\CP}{{\mathcal P}}  
\newcommand{\CQ}{{\mathcal Q}}  
\newcommand{\CR}{{\mathcal R}}
\newcommand{\CS}{{\mathcal S}}
\newcommand{\CT}{{\mathcal T}}
\newcommand{\CU}{{\mathcal U}}

\newcommand{\CX}{{\mathcal X}}

\newcommand{\fx}{{\mathfrak x}}
\newcommand{\fy}{{\mathfrak y}}

\newcommand{\FD}{{\mathfrak D}}

\newcommand{\FS}{{\mathfrak S}}

\newcommand{\BR}{{\mathbb R}}

\newcommand{\BC}{{\mathbb C}}

\newcommand{\BS}{{\mathbb S}}
\newcommand{\BZ}{{\mathbb Z}}
\newcommand{\BQ}{{\mathbb Q}}

\newcommand{\USL}{\CU_q(\fsl(2,\BR))}

\newcommand{\TP}{\tilde{\Phi}}
\newcommand{\TX}{\tilde{\Xi}}
\newcommand{\TI}{\tilde{\CI}}
\newcommand{\tf}{\tilde{f}}
\newcommand{\ux}{\underline{x}}
\newcommand{\uy}{\underline{y}}
\renewcommand{\Re}{\text{Re}}
\renewcommand{\Im}{\text{Im}}

\DeclareMathOperator*{\Res}{Res}

\def\ew{\hspace*{-1mm}}   \def\ppe{\hspace*{-2.5mm}}

\newcommand{\CBls}[6]{{#6}_{{\scriptstyle #1}}^s
  \hspace*{.3mm}\displaystyle{[} \ew \begin{array}{ll} {\scriptstyle #2 }
  \ppe & {\scriptstyle #3} \ppe \\[-2mm] {\scriptstyle #4}\ppe &
  {\scriptstyle #5}\ew \end{array}\displaystyle{]}}
\newcommand{\CBlt}[6]{{#6}_{{\scriptstyle #1}}^t
  \hspace*{.3mm}\displaystyle{[} \ew \begin{array}{ll} {\scriptstyle #2 }
  \ppe & {\scriptstyle #3} \ppe \\[-2mm] {\scriptstyle #4}\ppe &
  {\scriptstyle #5}\ew \end{array}\displaystyle{]}}
\newcommand{\CBlf}[6]{{#6}_{{\scriptstyle #1}}^{\flat}
  \hspace*{.3mm}\displaystyle{[} \ew \begin{array}{ll} {\scriptstyle #2 }
  \ppe & {\scriptstyle #3} \ppe \\[-2mm] {\scriptstyle #4}\ppe &
  {\scriptstyle #5}\ew \end{array}\displaystyle{]}}

\newcommand{\CpBl}[5]{\Phi_{{\scriptstyle #1}}^{\flat}
  \hspace*{.3mm}\displaystyle{[} \ew \begin{array}{ll} {\scriptstyle #2 }
  \ppe & {\scriptstyle #3} \ppe \\[-2mm] {\scriptstyle #4}\ppe &
  {\scriptstyle #5}\ew \end{array}\displaystyle{]}}
\newcommand{\CpBls}[5]{\Phi_{{\scriptstyle #1}}^s
  \hspace*{.3mm}\displaystyle{[} \ew \begin{array}{ll} {\scriptstyle #2 }
  \ppe & {\scriptstyle #3} \ppe \\[-2mm] {\scriptstyle #4}\ppe &
  {\scriptstyle #5}\ew \end{array}\displaystyle{]}}
\newcommand{\CpBlt}[5]{\Phi_{{\scriptstyle #1}}^t
  \hspace*{.3mm}\displaystyle{[} \ew \begin{array}{ll} {\scriptstyle #2 }
  \ppe & {\scriptstyle #3} \ppe \\[-2mm] {\scriptstyle #4}\ppe &
  {\scriptstyle #5}\ew \end{array}\displaystyle{]}}

\newcommand{\CGC}[6]{\displaystyle{[} \,\ew \begin{array}{lll} 
  {\scriptstyle #1} \ppe
  & {\scriptstyle #3} \ppe & {\scriptstyle #5} \ew \\[-2mm] {\scriptstyle
  #2} \ppe & {\scriptstyle #4}\ppe & {\scriptstyle #6} \ew\end{array}
  \displaystyle{]}}
\newcommand{\CGCZ}[6]{Z\displaystyle{(} \ew \begin{array}{lll} 
  {\scriptstyle #1} \ppe
  & {\scriptstyle #3} \ppe & {\scriptstyle #5} \ew \\[-2mm] {\scriptstyle
  #2} \ppe & {\scriptstyle #4}\ppe & {\scriptstyle #6} \ew\end{array}
  \!\displaystyle{)}}
\newcommand{\SJS}[6]{ \displaystyle{\bigl\{ } \ew 
\begin{array}{ll} {\scriptstyle #1 }
  \ppe & {\scriptstyle #2} \ppe \\[-2mm] {\scriptstyle #3}\ppe &
  {\scriptstyle #4}\ew \end{array}\big| \ew
\begin{array}{l} {\scriptstyle #5 }
  \ppe \\[-2mm] {\scriptstyle #6}\ew  \end{array}\displaystyle{\bigr\}_b}} 
 
\newcommand{\rf}[1]{(\ref{#1})}

\newcommand{\aufz}
{\begin{list}{$\bullet$}{\topsep0cm \itemsep0cm \parsep0cm}}
\newcommand{\eaufz}{\end{list}}
\newcounter{num}
\newcommand{\remlst}{\begin{list}
{(\arabic{num})}{\usecounter{num}\topsep0cm \itemsep0cm \parsep0cm}}
\begin{document}
\thispagestyle{empty}
\date{July 2000}
\title[Clebsch-Gordan and Racah-Wigner coefficients for representations of 
$\USL$]{Clebsch-Gordan and Racah-Wigner coefficients for a continuous
series of representations of $\USL$}
\author{\sc B. Ponsot, J. Teschner}
\address{B.P.: Laboratoire de Physique Math\'ematique,\\
Universit\'{e} Montpellier II,\\
Pl. E. Bataillon, 34095 Montpellier,\\ France}
\email{ponsot@lpm.univ-montp2.fr}
\address{J.T.: Institut f\"ur theoretische Physik,\\
Freie Universit\"at Berlin,\\
Arnimallee 14,\\ 
14195 Berlin, Germany}
\email{teschner@physik.fu-berlin.de}

\begin{abstract} The decomposition of tensor products of representations
into irreducibles is studied for a continuous family
of integrable operator representations
of $U_q(sl(2,\BR)$. It is described by an explicit
integral transformation involving a distributional kernel that can be seen
as an analogue of the Clebsch-Gordan coefficients. 
Moreover, we also study the relation between two canonical 
decompositions of triple tensor products into irreducibles. 
It can be represented by an integral transformation with a kernel 
that generalizes the Racah-Wigner coefficients. This kernel is 
explicitly calculated.
\end{abstract}

\maketitle
\section{Introduction}

Noncompact quantum groups can be expected 
to lead to very interesting generalizations 
of the rich and beautiful subject of
harmonic analyis on noncompact groups. 
Important progress has recently been made concerning 
an abstract ($C^*$-algebraic) theory of noncompact quantum
groups, see \cite{KV} for a nice overview and further references.
However, an important problem is still the rather limited 
supply of interesting examples. Results on the harmonic analysis
are so far only known for the quantum deformation of the group
of motions on the euclidean plane\cite{W1,DW},
the quantum Lorentz group \cite{BR1,BR2} and 
$SU_q(1,1)$ \cite{Kk}\cite{KSR}.
Moreover, 
there sometimes exist subtle analytical obstacles to 
construct quantum deformations of 
classical groups such as $SU(1,1)$ on the $C^*$-algebraic
level, cf. \cite{W2}.

Recently some evidence was presented in \cite{PT} that a certain noncompact 
quantum group with deformation parameter $q=e^{\pi i b^2}$ should 
describe a crucial internal structure of Liouville theory, 
a two-dimensional conformal field theory (CFT) that can be seen to be
as much a prototype for a CFT with continuous spectrum of 
Virasoro representations
as the harmonic analysis on
$SL(2,\BC)$ is a protoype for noncompact groups. 
The relation between Liouville theory and that quantum group
which was proposed in \cite{PT} generalizes the known 
equivalences between fusion categories of chiral algebras in
conformal field theories and braided tensor categories of quantum 
group representations, cf. e.g. \cite{KL,F}. These equivalences 
concern the isomorphisms that represent the operation of commuting tensor
factors as well as the associativity of tensor products, and can
be boiled down to the comparison of certain numerical data, 
the most non-trivial being some generalization of the Racah-Wigner 
coefficients (or fusion coefficients in CFT terminology).

The quantum group in question is $\USL$. A class of ``well-behaved''  
representation of $\USL$ on Hilbert-spaces was defined and classified
in \cite{S1}.
We will study a certain subclass of the representations
listed there. Some of the representations found in \cite{S1} 
reproduce known representions of principal or discrete series
of $\fsl(2,\BR)$ in the classical limit $b\ra 0$, others do not 
have a classical limit at all. The representations we will 
consider are of the latter type. Let us remark that representations
that are essentially equivalent to the class of representations
dicussed in our paper were recently also discussed in \cite{F1}.
The main result of the latter paper is a very interesting 
proposal for a braiding operation on such representations.

In our present paper we will present
explicit descriptions for the
decomposition of tensor products of these representations into
irreducibles, as well as the isomorphism relating two canonical
bases for triple tensor products. What appears to be remarkable
is the fact that the subseries we have picked out is 
actually closed under forming tensor products, which one 
would generally not expect if there exist other unitary 
representation. 
The maps describing the decomposition of tensor products lead to the
definition and explicit calculation of the generalization of the 
Racah-Wigner coefficients which represent the central ingredient
for the approach of \cite{PT} from the mathematics of quantum groups.

From the mathematical point of view one may view our results
as providing a technical basis for further studies of
a $C^*$ algebraic quantum group that may be generated\footnote{In
a similar sense as the bounded operators on $L^2(\BR)$ are 
generated by the unbounded operators $p$ and $q$ that satisfy
${[}p,q{]}=-i$, cf. \cite{W3} for more details} from
$\USL$ and its dual object, which is 
expected to be a $C^*$ algebraic quantum group generated
from $SL_q(2,\BR)$. In \cite{PT} we presented the definition of
$SL_q^+(2,\BR)$ as a quantum space, 
a $C^*$ algebra $\CA^+$ that is generated from $SL_q(2,\BR)$ and is 
acted on by analogues of left and right regular representation
of $\USL$. An $L^2$-space was introduced there, and the result
describing its decomposition into irreducible representations
of $\USL$ (Plancherel decomposition) was announced. 

Two aspects
of these constructions were unusual: $\CA^+$ was introduced such that
the elements $a,b,c,d$ generating $SL_q(2,\BR)$ have {\it positive}
spectrum and the $L^2$-space was introduced by a measure that has
no classical $q\ra 1 $ limit. It turns out that it is {\it precisely} the
subset of unitary $\USL$ representations studied in the present
paper which appears in the Plancherel decomposition of that $L^2$-space. 
We view these 
results as hints towards existence of a rather interesting 
$C^*$-algebraic quantum group related to $SL_q(2,\BR)$ that has no
classical counterpart, but other beautiful properties such as
a self-duality under $b\ra b^{-1}$ which are crucial for
the application to Liouville theory \cite{PT}. 

A first hint towards this self-duality can be
found in the observation made in \cite{PT}\cite{F1} (see
also \cite{F2} for closely related earlier observations) 
that the representations
that we consider may alternatively be seen as representations of 
$\CU_{\tilde{q}}(\fsl(2,\BR))$, where $\tilde{q}=e^{\pi i/b^2}$.
This led L. Faddeev to the proposal \cite{F1} to unify $\USL$ and 
$\CU_{\tilde{q}}(\fsl(2,\BR))$ into an object called 
``modular double'', which exhibits the self-duality under $b\ra b^{-1}$
in a manifest way. And indeed, it is found in the 
present paper that the Clebsch-Gordon intertwining maps, as well as the 
Racah-Wigner coefficients can be constructed in terms of a 
remarkable special
function $S_b(x)$. This special function
is closely related to the Barnes Double Gamma function \cite{Ba}, 
and was more recently independently introduced under the names
of ``Quantum Dilogarithm'' in \cite{FK}, and as
``Quantum Exponential function'' in \cite{Wo}. 
The function $S_b(x)$ has the property to be self-dual in the 
sense that it satisfies $S_b(x)=S_{1/b}(x)$. It follows 
from this self-duality of the function $S_b$ that 
the Clebsch-Gordan maps constructed in the present paper can be 
seen as intertwining maps for the ``modular double'' of L. Faddeev.

We would finally like to point out that our techniques for
dealing with finite difference operators that involve shifts
by imaginary amounts, in particular the method for 
determining the spectrum of such an operator, seem to be new 
and should have generalizations to a variety of other
problems where such operators appear. Moreover, the investigation
of the class of special functions that we use is fairly recent, 
so we will need to deduce several previously unknown properties.

The paper is organized as follows: In the following section 
we will introduce some technical preliminaries. Since we have to deal
with finite difference operators that shift the arguments of functions
by {\it imaginary} amounts, a lot of what follows will be based on the
theory of functions analytic in certain strips around the real axis, and 
the description of their Fourier-transforms via results of Paley-Wiener type.

The third section introduces the class of representations 
that will be studied in the present paper and discusses some of their
properties. 

This is followed by a section describing the decomposition of tensor
products of representations into irreducibles.

We then define and calculate b-Racah Wigner coefficients as the 
kernel that appears in the integral transformation that establishes
the isomorphism between two canonical decompositions of triple
tensor products.

Appendix A is in some sense the technical heart of the paper: It contains
the spectral analysis of a finite difference operator of second order
that is related to the Casimir on tensor products of two representations.
 
Appendices B and C contain some information on the special functions
that are used in the body of the paper.

{\bf Acknowledgements}
B.P. was supported in part by the EU under 
contract ERBFMRX CT960012. 
J.T. is supported by DFG SFB 288 ``Differentialgeometrie 
und Quantenphysik''.  Most of this work was carried out
while the second named author was 
at the Dublin Institute for Advanced Studies. He would like
to express this institution 
his sincere gratitude for support and hospitality.

\section{Preliminaries}

We collect some basic conventions, definitions and standard results that will 
be used throughout the paper.

\subsection{Finite difference operators}

The quantum group will be realized in terms of finite difference operators
that shift the arguments by an {\it imaginary} amount.
On functions $f(x)$, $x\in\BR$ 
that have an analytic continuation to a strip 
containing $\{x\in\BC;\Im(x)\in{[}-a_-,a_+{]}\}$, $a_{\pm}\geq 0$ one may 
define the finite difference operators $T_x^{ia}$, $a\in{[}-a_-,a_+{]}$ by
\begin{equation}
T_x^{ia}f(x)=f(x+ia).
\end{equation}
As convenient notation we will use
\begin{equation}
[x]_b^{}\equiv
\frac{\sin(\pi b x)}{\sin(\pi b^2)}, \qquad d_x\equiv\frac{1}{2\pi}\pa_x,
\qquad 
[d_x+a]_b^{}\equiv\frac{e^{\pi i b a}T_x^{\frac{ib}{2}}-
e^{-\pi i b a}T_x^{-\frac{ib}{2}}}{e^{\pi i b^2}-e^{-\pi i b^2}}.
\end{equation}

\subsection{Fourier-transformation}

Our
notation and conventions concerning the Fourier-transformations are as 
follows:
Let $\CS(\BR)$ denote the usual Schwartz-space of functions on the real
line. The Fourier-transformation of a function $f\in\CS(\BR)$ will be
defined as
\begin{equation}
\tilde{f}(\om)=\int\limits_{-\infty}^{\infty}dx\;e^{-2\pi i \om x}f(x).
\end{equation}
The corresponding inversion formula is then
\begin{equation}\label{FT}
f(x)=\int\limits_{-\infty}^{\infty}d\om\;e^{2\pi i \om x}\tilde{f}(\om).
\end{equation}

The Fourier-transformation maps the finite difference operator
$T_x^{ia}$ to the operator of multiplication with $e^{-2\pi a\om}$. It will
therefore be a useful tool for dealing with these operators. 
Of fundamental importance will be the connection between   
 analyticity of functions
in a strip to exponential decay properties of its Fourier-transform
and vice versa that is expressed by the classical Paley-Wiener theorem:

\begin{thm} {(Paley-Wiener)}
Let $f$ be in $L^2(\BR)$. Then $(e^{2\pi xa_+}+e^{-2\pi xa_-})f\in L^2(\BR)$,
$a_{\pm}>0$ if and only if $\tilde{f}$ has an analytic continuation to the 
strip $\{\om\in\BC;\Im(\om)\in(-a_-,a_+)\}$ such that for any 
$\om_2\in(-a_-,a_+)$, $\tilde{f}(.+i\om_2)\in L^2(\BR)$ and 
\begin{equation}
\sup_{\om_2\leq b}\;\,\int\limits_{-\infty}^{\infty}d\om_1\;
|\tilde{f}(\om_1+i\om_2)|^2 \;\; <\;\;\infty \quad\text{for any} \quad
b\in(-a_-,a_+).\end{equation}
\end{thm}
\begin{proof}Cf. e.g. \cite{K}.\end{proof}
The following simple variant of this result will often be useful:
\begin{lem}
\label{PW}
For $f\in \CS(\BR)$, the following two conditions are equivalent:
\begin{enumerate}
\item $f$ 
is the restriction to $\BR$ of a function $F$ that is meromorphic 
in the strip $\{ z\in\BC;\Im(z)\in(-a_-,a_+)\}$, $a_+,a_->0$ with 
finitely many poles
in the upper (lower) half plane at $\CP_{\pm}\equiv \{z_j;j\in\CI_{\pm}\}$,
$|{\rm Im}(z_j)|>0$, and 
all functions $F_y(x)\equiv F(x+iy)$, $y\in(-a_-,a_+)$ are of
rapid decrease, and
\item one has the following asymptotic behavior 
of the Fourier-transform $\tilde{f}(\om)$ for $\om\ra\pm \infty$:
\begin{equation*}
\begin{aligned}
\tilde{f}(\om)\;=\;&  -2\pi i \sum_{j\in\CI_{-}}e^{-2\pi iz_j\om}\Res_{z=z_j}
F(z)+\tilde{f}_{a_+}(\om) \\
\tilde{f}(\om)\;=\;&  +2\pi i \sum_{j\in\CI_{+}}e^{-2\pi iz_j\om}\Res_{z=z_j}
F(z)+\tilde{f}_{a_-}(\om),
\end{aligned}
\end{equation*}
where $\tilde{f}_{a_{\pm}}(\om)$ 
decay as $x\ra\pm\infty$ faster than $e^{-2\pi a|\om|}$ for 
any $a\in(-a_-,a_+)$.
\end{enumerate}
\end{lem}
\subsection{Distributions}

Let $\CS'(\BR)$ be  
the space of tempered distributions on $\CS(\BR)$. 
The dual pairing between a distributions $\Phi\in\CS'(\BR)$
and a function $f\in\CS(\BR)$ will be denoted by $\bra \Phi,f\ket$.
The Fourier transformation on $\CS'(\BR)$ is defined by 
$\bra \tilde{\Phi},\tilde{f}\ket
\equiv \bra \Phi,f\ket$ for any 
$f\in \CS(\BR)$. It should be noted that if a distribution $\Phi\in\CS'(\BR)$
actually happens to be represented by a function $\Phi(x)$ via
\[
\bra \Phi,f\ket = \int\limits_{-\infty}^{\infty}dx\; \Phi(x)f(x)
\]
then our definition of the Fourier-transform of $\Phi$ implies that
instead of \rf{FT} one has the following inversion formula for $\Phi(x)$: 
\begin{equation}\label{distFT}
\Phi(x)=\int\limits_{-\infty}^{\infty}d\om\;e^{-2\pi i \om x}\tilde{\Phi}(\om).
\end{equation}

The distributions that appear below will all be defined in terms of 
meromorphic functions by means of the so-called $i\ep$-prescription:
Assume given a familiy of functions $\Phi_{\ep}$, $\ep>0$ that are 
meromorphic in some strip containing $\BR$, rapidly
decreasing at infinity and have finitely many poles
with $\ep$-independent residues at a distance $\ep$ from the real axis. 
The limit $\Phi\equiv \lim_{\ep\ra 0}
\Phi_{\ep}$ then defines a distribution $\Phi\in \CS'(\BR)$. We
will often use the symbolic notation $\Phi(x)$ for the resulting
distribution, keeping in mind that $\Phi(x)$ will not be defined for
all $x\in\BR$.

There is a simple generalization of Lemma \ref{PW} to 
such distributions in $\CS'(\BR)$: Poles on the real axis
correspond to asymptotic behavior of the form $e^{2\pi i \om x}$
of the Fourier-transform: 
\begin{lem}\label{distPW} 
For $\Phi\in \CS'(\BR)$, the following two conditions are equivalent:
\begin{enumerate}
\item $\Phi=\lim_{\ep\ra 0}\Phi_{\ep}$, where $\Phi_{\ep}$ is for 
$\ep>0$ represented as the restriction to $\BR$ of a function 
$\Phi_{\ep}(x)$ that is meromorphic 
in the strip $\{ z\in\BC;\Im(z)\in(-a_-,a_+)\}$, $a_+,a_->0$ with 
finitely many poles
in the upper (lower) half plane at $\CP_{\pm}^{\ep}\equiv 
\{z_j\pm i\ep ;j\in\CI_{\pm}\}$,
$\pm{\rm Im}(z_j)\geq 0$, and 
all functions $\Phi_{\ep,y}(x)\equiv \Phi_{\ep}(x+iy)$, $x,y\in\BR$,
$y\in(-a_+,a_-)$ are of
rapid decrease, and
\item $\TP$ is represented by a function $\TP(\om)\in\CC^{\infty}(\BR)$ that 
has the following asymptotic behavior:
\begin{equation*}
\begin{aligned}
\tilde{\Phi}(\om)\;=\;&  +2\pi i \sum_{j\in\CI_{+}}e^{2\pi iz_j\om}\Res_{z=z_j}
\Phi(z)+\tilde{\Phi}_{a_+}(\om) \\
\tilde{\Phi}(\om)\;=\;&  -2\pi i \sum_{j\in\CI_{-}}e^{2\pi iz_j\om}\Res_{z=z_j}
\Phi(z)+\tilde{\Phi}_{a_-}(\om),
\end{aligned}
\end{equation*}
where $\TP_{a_{\pm}}(\om)$ decay faster than than $e^{-2\pi a|\om|}$ for 
any $a\in(-a_-,a_+)$.
\end{enumerate}
\end{lem} 
\begin{rem}
The sign flips between Lemma \ref{PW} and Lemma 
\ref{distPW} are due to the different inversion formulae for functions
and distributions.
\end{rem}

\subsection{A useful Lemma from complex analysis}

The following Lemma is useful for determining the analytic properties
of convolutions of meromorphic functions:
\begin{lem} \label{doubpole}
Let $f(z_0;z_1,z_2)$ be meromorphic in its variables
in some open strip $\CS$ around the real axis, with singular 
behavior near $z_0=z_1=z_2$ of the form $R_{12}(z_1)
(z_0-z_1)^{-1}(z_0-z_2)^{-1}$.
The function $I(z_1,z_2)$, defined by the integral 
\begin{equation}\label{Idef}
I(z_1,z_2)\;\equiv\;
\int\limits_{-\infty}^{\infty}dz_0
\;\,f(z_0;z_1,z_2),
\end{equation}
will then be a function that
has a meromorphic continuation w.r.t. $z_i$, $i=1,2$ 
to the whole strip $\CS$. If $z_1$ and $z_2$
were initially seperated by the real axis one will find 
a pole with residue $R_{12}(z_1)$
at $z_1=z_2$. If not, $I(z_1,z_2)$ will be nonsingular at $z_1=z_2$
as well.
\end{lem}
\begin{proof}
To define the meromorphic continuation of $I(z_1,z_2)$ in cases 
where the poles $z_i$, $i=1,2$ cross the contour of integration 
of the integral \rf{Idef} one just needs to deform the contour
accordingly. This will obviously always be possible as long as 
$z_i$, $i=1,2$ were initially not separated by the real axis.
We will therefore
turn to the case that they were initially seperated, and consider
w.l.o.g. the case that $z_1$ was initially in the upper,
$z_2$ in the lower half plane. In this case one may deform the 
contour into a contour that passes {\it above} $z_1$ plus a small 
circle around $z_1$. The residue contribution from the integral
over that small circle is
\begin{equation}
2\pi i \frac{R_{12}(z_1)}{z_1-z_2}+({\rm contributions\;\,
regular\;\, as}\;\, z_1-z_2\ra 0)
\end{equation}
The Lemma is proven. \end{proof}

\section{A class of representations of $U_q(sl(2,\BR))$}

\subsection{Definintion}

$U_q(sl(2,\BR)$ 
is a Hopf-algebra with 
\begin{equation}
\begin{aligned}
{}& \text{generators:}\quad  E,\quad F,\quad K,\quad K^{-1};\\
& \text{relations:}\quad 
KE=qEK,\qquad\quad KF=q^{-1}FK,\qquad[E,F]=-\frac{K^2-K^{-2}}{q-q^{-1}};\\
& \text{star-structure:}\quad 
K^*=K, \qquad E^*=E, \qquad F^*=F;\\ 
& \text{co-product:}\quad\De(K)=K\ot K,\qquad 
\begin{aligned}\De(E)=&E\ot K+K^{-1}\ot E,\\
\De(F)=&F\ot K+K^{-1}\ot F.
\end{aligned}\end{aligned}\end{equation}
The center of $U_q(sl(2,\BR)$ is generated by the $q$-Casimir
\begin{equation}
C=FE-\frac{qK^2+q^{-1}K^{-2}-2}{(q-q^{-1})^2}.
\end{equation}
We will consider the case that $q=e^{\pi i b^2}$, 
$b\in(0,1)\cap(\BR\setminus\BQ)$.

Unitary representations of $U_q(sl(2,\BR))$ by operators on a Hilbert-space
have been studied in \cite{S1}. Since there are no unitary representations
in terms of bounded operators some care is needed
in order to single out an interesting 
class of ``well-behaved'' representations. A natural notion
of ``well-behaved'' was introduced in \cite{S1}, where the 
corresponding unitary representations of $U_q(sl(2,\BR))$ were 
classified. 

In the present paper we will study a one-parameter subclass $\CP_{\al}$,
$\al\in Q/2+i\BR$, $Q=b+b^{-1}$
of the representations listed in \cite{S1} which are constructed as 
follows: The representation will be realized on the
space $\CP_{\al}$ of entire analytic functions $f(x)$ that 
have a Fourier-transform $f(\om)$ which is meromorphic in $\BC$
with possible poles at
\begin{equation}\begin{aligned}
\om\;=\;&i(\al-Q-nb-mb^{-1})\\
\om\;=\;&i(Q-\al+nb+mb^{-1})
\end{aligned}\qquad n,m\in\BZ^{\geq 0}.\end{equation}

\begin{rem} It can be shown that $\CP_{\al}$ is a Frechet-space.
\end{rem}

One may then introduce the following finite difference operators 
\begin{equation}\label{uslgens}
\begin{aligned}
\pi_{\al}(E)\;\equiv\;&e^{+2\pi b x}
[d_x+Q-\al]_b^{}\\
\pi_{\al}(F)\;\equiv\;&e^{-2\pi b x}
[d_x+\al-Q]_b^{}
\end{aligned}\qquad\qquad \pi_{\al}(K)\;\equiv\;T_x^{\frac{ib}{2}}.
\end{equation}
As shorthand notation we will also use $u_{\al}\equiv\pi_{\al}(u)$.
\begin{lem}

\begin{itemize}
\item[(i)] The operators $\pi_{\al}(u)$, $u=E,F,K$ map 
$\CP_{\al}$ into itself.
\item[(ii)] $\pi_{\al}(u)$, $u=E,F,K$ generate 
a representation of $\USL$ on $\CP_{\al}$.
\end{itemize}\end{lem}
\begin{proof} 
To verify (i), note that Fourier-transformation maps $E_{\al}$,
$F_{\al}$, $K_{\al}$ into the following operators:
\begin{equation}\label{uslgenft}
\begin{aligned}
\tilde{E}_{\al}=&
[-i\om+\al]_b^{} T_{\om}^{ib}\\
\tilde{F}_{\al}=& 
[-i\om-\al]_b^{} T_{\om}^{-ib}
\end{aligned}\qquad\qquad K_{\al}=e^{-\pi b\om}.
\end{equation}
The claim follows from the fact that $[x]_b=0$ for $x=nb^{-1}$, $n\in\BZ$.

(ii) is checked by straightforward calculation.
\end{proof}
\begin{propn}\label{intreps}
The operators \rf{uslgens} generate an {\rm integrable} operator representation
of $\USL$ in the sense of \cite{S1}, i.e. \begin{enumerate}
\item $E_{\al}$, $F_{\al}$, $K_{\al}$ have self-adjoint extensions in
$L^2(\BR)$,
\item the corresponding unitary operators $E^{it}_{\al}$, $F^{it}_{\al}$, 
$K^{it}_{\al}$
satisfy 
\[ K^{is}_{\al}E^{it}_{\al}=q^{-ts}E^{it}_{\al}K^{is}_{\al},\qquad
K^{is}_{\al}F^{it}_{\al}=q^{ts}F^{it}_{\al}K^{is}_{\al},\quad\text{and} \]
\item the q-Casimir strongly commutes with $E_{\al}$, $F_{\al}$ and $K_{\al}$.
\end{enumerate}\end{propn}
\begin{proof}
It suffices to show that the representation $\CP_{\al}$ is unitarily 
equivalent to one of the representations listed in \cite{S1}.
Consider the operator $J_{\al}$ defined as 
$(J_{\al}\tilde{f})(\om)=S_b(\al-i\om)\tilde{f}(\om)$ in terms
of the special function $S_b(x)$ (cf. Appendix B). 
$J_{\al}$ is unitary since $|S_b^{-1}(\al-i\om)|^2=1$ which follows 
from eqn. \rf{reflprop} in Appendix B. Moreover, it follows from the 
analytic and asymptotic properties of $S_b(x)$ given in the 
Appendix that $J_{\al}$ maps $\CP_{\al}$ to the space 
$\CR_{\al}$ of entire analytic functions which have a Fourier-transform
that is meromorphic in $\BC$ with possible poles at
\begin{equation}\begin{aligned}
\om\;=\;&i(\al-Q-nb-mb^{-1})\\
\om\;=\;&i(-\al-nb-mb^{-1})
\end{aligned}\qquad n,m\in\BZ^{\geq 0}.\end{equation} 
One finally finds from the functional relations of the $S_b$-functions,
eqn. \rf{S:funrel} that
\begin{equation}
\begin{aligned}
J_{\al}^{-1}\tilde{E}_{\al}^{}J_{\al}^{}=&
T_{\om}^{ib}\\
J_{\al}^{-1}\tilde{F}_{\al}^{}J_{\al}^{}=& 
[\al+i\om]_b^{}T_{\om}^{-ib}[\al-i\om]_b^{} 
\end{aligned}\qquad\qquad J_{\al}^{-1}K_{\al}^{}J_{\al}^{}=e^{-\pi b\om}.
\end{equation}
Our representation is thereby easily recognized as 
the representation denoted by $(I)_{1,-1,c}$ in Corollary 5 of \cite{S1}, 
where
$c=[\al-\frac{Q}{2}]_b^2+2(q-q^{-1})^{-2}$. Note that
our notation $Q$ is different from that in \cite{S1} and
$c\leq 2(q-q^{-1})^{-2}$.
\end{proof}

\begin{rem}
The representations considered here 
form a subset of the representations
of $\CU_q(\fsl(2,\BR))$ that appear in the classification of \cite{S1}. 
This subset has the following remarkable property:
If one introduces generators $\tE$, $\tF$, $\tK$ by 
replacing $b\ra b^{-1}$ in the
expressions for $E$, $F$, $K$ given above, one obtains 
a representation of $\CU_{\tilde{q}}(\fsl(2,\BR))$
$\tilde{q}=\exp(\pi i b^{-2})$ on the same space $\CP_{\al}$.
The generators $\tE$, $\tF$, $\tK^2$ commute with
$E$, $F$, $K^2$ on the space $\CP_{\al}$. This does not mean, 
however, that these operators commute as self-adjoint operators
on $L^2(\BR)$. This self-duality property of our representations $\CP_{\al}$
is related to the fact that the representations 
$(\CP_{\al},\pi_{\al})$ do {\it not} have a classical ($b\ra 0$) 
limit.
\end{rem}

\subsection{Intertwining operators}

The representations with labels $\al$ and $Q-\al$ are equivalent.
The unitary operator establishing this equivalence can be most easily
found by considering the Fourier-transform of the representation \rf{uslgens},
as already done in the proof of Proposition \ref{intreps}, eqns. \rf{uslgenft}:
Define the operator $\TI_{\al}:L^2(\BR)\ra L^2(\BR)$ as 
\begin{equation}
(\TI_{\al}\tf)(\om)\;=\;\TB_{\al}(\om)f(\om),\quad \TB_{\al}(\om)
\;\equiv\; \frac{S_b(\al-i\om)}{S_b(Q-\al-i\om)} .
\end{equation}
The operator $\TI_{\al}$ is unitary since $|\TB_{\al}(\om)|=1$. It maps 
$\CP_{\al}$ to $\CP_{Q-\al}$ as follows
from the analytic and asymptotic
properties of the $S_b$-function summarized in Appendix B.
The fact that
\begin{equation}
\pi_{Q-\al}(u)\TI_{\al}=\TI_{\al}\pi_{\al}(u),\quad
u\in\USL
\end{equation}
is a simple consequence of the functional relations \rf{S:funrel}, Appendix
B of the $S_b$-functions. 

By inverse Fourier-transformation one finds the representation of 
the intertwining operator on functions $f(x)$. It takes the form
\begin{equation}
(\CI_{\al}f)(x)\;=\;\int\limits_{\BR}dx'\;B_{\al}(x-x')f(x),
\end{equation}
where the inverse Fourier-transform defining the kernel $B_{\al}(x-x')$
may be found by means of eqn. \rf{bbeta}, Appendix B to be given by
\begin{equation}
B_{\al}(x-x')\;=\;S_b(2\al)\frac{S_b\bigl(\frac{Q}{2}+i(x-x')-\al\bigl)}
{S_b\bigl(\frac{Q}{2}+i(x-x')+\al\bigl)}.
\end{equation}

\section{The Clebsch-Gordan decomposition of tensor products}

The co-product allows us to define the tensor product of representations:
For any $u\in\USL$ let $\pi_{21}(u)\equiv (\pi_{\al_2}\ot\pi_{\al_1})\De(u)$.
The operators $\pi_{21}(u)$ generate a representation of $\USL$
on $\CP_{\al_2}\ot\CP_{\al_1}$. Our aim is to determine the 
decomposition of this representation into irreducible 
representations of $\USL$.
\begin{lem} $\CP_{\al_2}\ot\CP_{\al_1}$ is 
dense in $L^2(\BR)\ot L^2(\BR)$.
\end{lem}
\begin{proof}
Any two-variable Hermite-function is contained in $\CP_{\al_2}\ot\CP_{\al_1}$.
\end{proof}  
\begin{defn}
Define a distributional kernel 
$\CGC{\al_3}{x_3}{\al_2}{x_2}{\al_1}{x_1}$
(the ``Clebsch-Gordan coefficients'') by an expression of the form
\begin{equation}
\CGC{\al_3}{x_3}{\al_2}{x_2}{\al_1}{x_1}\;\equiv
\; \lim_{\ep\downarrow 0}\; \CGC{\al_3}{x_3}{\al_2}{x_2}{\al_1}{x_1}_{\ep}^{},
\end{equation}
where the meromorphic function
$\CGC{\al_3}{x_3}{\al_2}{x_2}{\al_1}{x_1}_{\ep}^{}$ is defined as
\begin{equation}\begin{aligned}\label{Clebsch}
\CGC{Q-\al_3}{x_3}{\al_2}{x_2}{\al_1}{x_1}_{\ep}
\;= \; & e^{-\frac{\pi i}{2}(\De_{\al_3}-\De_{\al_2}-\De_{\al_1})}\\
& \qquad\ti
D_b(\be_{32};y_{32}+\ep)D_b(\be_{31};y_{31}+\ep)D_b(\be_{21};y_{21}+\ep), 
\end{aligned}\end{equation}
$\De_{\al}=\al(Q-\al)$, the distribution $D_b(\al;y)$ is 
defined in terms of the 
Double Sine function $S_b(y)$ (cf. Appendix) as 
\begin{equation}
D_b(\al;y)=\frac{S_b(y)}{S_b(y+\al)},
\end{equation}
and the coefficients $y_{ji}$, $\be_{ji}$, $j>i\in\{1,2,3\}$ are given by
\begin{equation} \label{CGcoeff}
\begin{aligned}
y_{32}=& i(x_3-x_2)-\fr{1}{2}(\al_3+\al_2-Q)\\
y_{31}=& i(x_1-x_3)-\fr{1}{2}(\al_3+\al_1-Q)\\
y_{21}=& i(x_1-x_2)-\fr{1}{2}(\al_2+\al_1-2\al_3)
\end{aligned}
\qquad \quad \begin{aligned}
\be_{32}=& \al_2+\al_3-\al_1\\
\be_{31}=& \al_3+\al_1-\al_2\\
\be_{21}=& \al_2+\al_1-\al_3.
\end{aligned}
\end{equation}
\end{defn}
The aim of this section will be to prove
\begin{thm}\label{Cldeco}
The $\CU_{q}(\fsl(2,\BR))$-representation 
$\pi_{21}$ defined on $\pi_{\al_2}\ot\pi_{\al_1}$ 
decomposes as follows into irreducible representations $\CP_{\al}$: 
\begin{equation}\label{CGdecoQG}
\pi_{\al_2}\ot \pi_{\al_1}\;\,\simeq\;\, \int\limits_{\BS}^{\oplus}
\!d\al\;\, \pi_{\al,}\qquad \BS\equiv \frac{Q}{2}+i\BR^+.
\end{equation}
The isomorphism can be described explicitly in terms of a 
unitary map $\CC_{21}$ of the form
\begin{equation}\label{CGmap}
\CC_{21}\quad : \quad  \begin{aligned}
L^2(\BR\times\BR) \quad &  \ra \quad
L^2(\BS\times\BR,d\mu(\al_3)dx_3),\qquad d\mu(\al)\equiv |S_b(2\al)|^2\\
f(x_2,x_1) \quad &  \ra \quad F_f(\al_3,x_3)\equiv 
\int_{\BR}dx_2dx_1\;\CGC{\al_3}{x_3}{\al_2}{x_2}{\al_1}{x_1}\;
f(x_2,x_1)
\end{aligned}
\end{equation}
such that the corresponding projections
$\Pi_{21}(\al_3)$, $\bigl(\Pi_{21}(\al_3)f\bigr)(x_3)=F_f(\al_3,x_3)$, map 
$\CP_{\al_2}\ot\CP_{\al_1}$ into $\CP_{\al_3}$ and
intertwine the respective $\USL$ actions according to
\begin{equation}\label{intertw}
\Pi_{21}(\al_3)\pi_{21}(u)=\pi_{\al_3}(u)
\Pi_{21}(\al_3)\qquad u\in\USL.
\end{equation}
\end{thm}
\begin{rem}
It follows from Theorem \ref{Cldeco} that the representation 
$\pi_{21}$ is in fact integrable, which was not clear
apriori.
\end{rem}
\begin{rem}
It is remarkable and nontrivial that the subset of ``self-dual'' integrable
representations of $\USL$ is actually closed under tensor products. 
\end{rem}
\begin{rem}
The appearance of the measure $d\mu(\al)$ is natural since
$d\mu(\al)$ is the Plancherel measure for the dual space of
functions $L^2(SL^+_q(2,\BR))$, cf. \cite{T}.
\end{rem}
\begin{cor}
The Clebsch-Gordan coefficients $\CGC{\al_3}{x_3}{\al_2}{x_2}{\al_1}{x_1}$
satisfy the following orthogonality and completeness relations:
\begin{equation}\begin{aligned}\label{CGorth}
{} & \lim_{\ep\downarrow 0}\;\int\limits_{\BR}dx_1dx_2\;\,
\CGC{\al_3}{x_3}{\al_2}{x_2}{\al_1}{x_1}^*_{\ep}
\CGC{\be_3}{y_3}{\al_2}{x_2}{\al_1}{x_1}^{}_{\ep}
=|S_b(2\al_3)|^{-2}\de(\al_3-\be_3)\de(x_3-y_3)\\
{} & \lim_{\ep\downarrow 0}\;\int\limits_{\BS}d\al_3 \; |S_b(2\al_3)|^2
\int\limits_{\BR}dx_3 \;\,
\CGC{\al_3}{x_3}{\al_2}{x_2}{\al_1}{x_1}^*_{\ep}
\CGC{\al_3}{x_3}{\al_2}{y_2}{\al_1}{y_1}_{\ep}^{}
=\de(x_2-y_2)\de(x_1-y_1). 
\end{aligned}\end{equation}
\end{cor}

The main step in the proof of Theorem \ref{Cldeco} will be the construction
of a common spectral decomposition for the operators 
$Q_{21}\equiv(\pi_{\al_2}\ot\pi_{\al_1})\De(Q)$ and $K_{21}$.
The decomposition of $L^2(\BR\times\BR)$ into eigenspaces of $K_{21}$
is simply obtained by Fourier-transformation:
\begin{equation}
\CF\quad : \quad  \begin{aligned}
L^2(\BR\times\BR) \quad &  \ra \quad
L^2(\BR\times\BR)\\
f(x_2,x_1) \quad &  \ra \quad F(\ka_3,x_-)\equiv 
\int_{\BR}dx_+\;\, e^{-\pi i \ka_3x_+}
f\bigl(\fr{x_++x_-}{2},\fr{x_+-x_-}{2}\bigr)
\end{aligned}
\end{equation} 
The q-Casimir $Q_{21}$ is mapped under this Fourier-transformation $\CF$ into 
a second order finite difference operator $C_{21}(\ka_3)$
that contains shifts w.r.t. 
the variable $x_-$ only and therefore leaves the eigenspaces of $K_{21}$
invariant: 
\begin{equation}\begin{aligned}
C_{21}(\ka_3) & -\bigl[\al_3-\fr{Q}{2}\bigr]_b^2=\\ =\quad
& {[}-ix-\fr{1}{2}(\al_1+\al_2-Q)+(\al_3-\fr{Q}{2}){]}_b^{}
                 {[}-ix-\fr{1}{2}(\al_1+\al_2-Q)-(\al_3-\fr{Q}{2}){]}_b^{}\\
-& {[}-ix+\fr{1}{2}(\al_1+\al_2)-Q{]}_b^{}\Bigl(
e^{i\pi b(-ix-\frac{1}{2}(\al_1+\al_2))}\{\al_1-\al_2+i\ka_3\}_b^{}\\
& \qquad\qquad\qquad\qquad\qquad\quad
-e^{-i\pi b(-ix-\frac{1}{2}(\al_1+\al_2))}\{\al_1-\al_2-i\ka_3\}_b\Bigl)
T_{x_-}^{-ib}\\
+& {[}-ix+\fr{1}{2}(\al_1+\al_2)-Q{]}_b^{}
{[}-ix+\fr{1}{2}(\al_1+\al_2)-2Q{]}_b^{} T_{x_-}^{-2ib},
\end{aligned}
\end{equation}
where the following notation has been used: 
\begin{equation}
[x]_b\equiv \frac{\sin(\pi b x)}{\sin(\pi b^2)}, \qquad
\{ x\}_b \equiv \frac{\cos(\pi b x)}{i\sin(\pi b^2)}.
\end{equation}
The spectral analysis of the operator $C_{21}$ is performed
in Appendix A. The result may be summarized as follows:
Eigenfunctions $\Phi_{\al_3}(\al_2,\al_1|\ka_3|x)$ of $C_{21}$ 
are given by an expression of the form
\begin{equation}\label{Phidef}
\Phi_{Q-\al_3}   (\al_2,\al_1|\ka_3|x)\;=   
 M^{\al_3;\ka_3}_{\al_2,\al_1}
\;e^{\pi x(2\al_3-2\al_2+i\ka_3)}\;\Theta_b(T,y_-)
\;\Psi_b(U,V,W;y_+).\end{equation}
The special functions $\Theta_b(T;y)$ and $\Psi_b(U,V,W;y)$ are defined
in Appendix B, $y_{\pm}$ are introduced as 
$y_{\pm}=-ix-\frac{1}{2}(\al_2+\al_1-Q)\mp(\al_3-\fr{Q}{2})$ 
and the coefficients $T$, $U$, $V$, $W$ are given as
\begin{equation}\begin{aligned}
T=& \al_2+\al_1-\al_3\\
U=& \al_3+\al_1-\al_2
\end{aligned}\qquad\begin{aligned} 
V=& -i\ka_3+\al_3\\
W=& -i\ka_3+\al_1-\al_2+Q.
\end{aligned} 
\end{equation}

\begin{thm}\label{casspec}
A complete set of generalized eigenfunctions for the operator $C_{21}(\ka_3)$
is given by $\{(\Phi_{\al_3})^* ;
\al_3\in\BS\}$.
\end{thm}

By combining Theorem \ref{casspec} with the usual Plancherel formula for the 
Fourier-transformation $\CF$ one concludes that each function 
$f(x_2,x_1)\in L^2(\BR\ti\BR)$ can be decomposed as 
($x_{\pm}\equiv x_2\pm x_1$)
\begin{equation}
f (x_2,x_1)= \int\limits_{\BR}d\ka_3 \;\,e^{\pi i \ka_3x_+}
\int\limits_{\BS}d\mu(\al_3) 
\;\,\bigl(\Phi_{\al_3}
(\al_2,\al_1|\ka_3|x_-)\bigr)^* F_f(\al_3,\ka_3),
\end{equation}
where the generalized Fourier-transformation $F_f$ of $f$ 
is defined as
\begin{equation}
F_f(\al_3,\ka_3)=\int\limits_{\BR}dx_2dx_1\;\,
e^{-\pi i \ka_3x_+}\;\Phi_{\al_3}
(\al_2,\al_1|\ka_3|x_-)f(x_2,x_1).
\end{equation}
The measure $d\mu(\al_3)$ will be determined later.
One may next observe that 
\begin{lem} One has 
\begin{equation}\label{FTk3}
\CGC{\al_3}{\ka_3}{\al_2}{x_2}{\al_1}{x_1}
\;\equiv \;\int\limits_{\BR}dx_3  \;\,e^{2\pi i \ka_3x_3}\; 
\CGC{\al_3}{x_3}{\al_2}{x_2}{\al_1}{x_1}\;=\;
e^{-\pi i \ka_3x_+}\;\Phi_{\al_3}
(\al_2,\al_1|\ka_3|x_-),
\end{equation}
if the normalization  factor $M$ in \rf{Phidef} is chosen as
\begin{equation}
M^{\al_3;\ka_3}_{\al_2,\al_1}
\;\equiv  \; e^{\pi i \al_2(\al_2-\al_3)}e^{-\pi i (\al_3-i\ka_3)(\al_3+\al_2
-Q)}
\end{equation}
\end{lem}
\begin{proof} 
The kernel $\CGC{Q-\al_3}{x_3}{\al_2}{x_2}{\al_1}{x_1}$ may be rewritten
in terms of the function $\Theta_b(\be;y)$ as
follows:
\begin{equation} \begin{aligned}
\CGC{Q-\al_3}{x_3}{\al_2}{x_2}{\al_1}{x_1}=
e^{\pi i \al_1\al_2}& 
e^{2\pi(x_3(\al_2-\al_1)+\al_1x_1-\al_2x_2)}\\
& \times
\Theta_b(\be_{32};y_{32})\Theta_b(\be_{31};y_{31})\Theta_b(\be_{21};y_{21}). 
\end{aligned}\end{equation}
The substitution $s=-i(x_3-x_2)+\frac{1}{2}(\al_3+\al_2-Q)$ then leads
to the Euler-type integral \rf{Eulerint} for the b-hypergeometric function.
The rest is straightforward.
\end{proof}
If follows that the generalized Fourier-transformation defined
in Theorem \ref{casspec} represents a decomposition into eigenspaces of the 
q-Casimir $Q_{21}$. Two things remain to be done in order to finish
the proof of Theorem \ref{Cldeco}: On the one hand it remains to calculate 
the spectral measure $d\mu(\al_3)$, and on the other hand one needs 
to verify the intertwining 
property \rf{intertw}.

\subsection{Spectral measure}

We will show in this subsection that $d\mu(\al_3)=|S_b(2\al_3)|^2$.
This follows from the combination of the following two results.
We first of all determine the asymptotics of 
the distributional
Fourier-transform of $\Phi_{\al_3}$:
\begin{lem}
The function $\TP_{\al_3}(\om)$ (defined as in \rf{distFT}) 
decays exponentially for $\om\ra\infty$ and has the following
asymptotic behavior for $\om\ra -\infty$:
\begin{equation}\label{TPasym}
\TP_{\al_3}(\om)\;=\;N_+(\al_3)e^{2\pi i \om x_+}+N_-(\al_3)
e^{2\pi i \om x_-}+R_-(\om),
\end{equation}
where $R_-(\om)$ decays exponentially for $\om\ra-\infty$, $x_+$ and $x_-$ 
are defined by
\[ x_{\pm}\equiv 
+\fr{i}{2}\bigl(\al_1+\al_2-Q\bigr)\pm i\bigl(\al_3-\fr{Q}{2}\bigr)
\]
and $|N_{\pm}(\al_3)|^2= |S_b(2\al_3)|^{-2}$. 
\end{lem}
\begin{proof}
According to Lemma \ref{distPW} one just needs to calculate
the residues of $\Phi_{\al_3}$ for the poles at $x=x_{\pm}$.
We will only need the absolute values of these quantities.

The pole at $x=x_-$ comes from the $G_b/G_b$ factor in the
expression for $\Phi$. To calculate its residue one needs
the following special value of the $\Psi$-function:
\begin{equation}
\Psi_b(U,V;W;W-U-V)=\frac{G_b(V)G_b(W-U-V)}{G_b(W-U)},
\end{equation} 
which follows easily from the fact that the representation \rf{Eulerint}
simplifies to the b-beta integral \rf{bbeta} for $x=W-U-W$. 
We furthermore note that $|G_b(\frac{Q}{2}+ix)|^2=1$ from the 
reflection property of $S_b(x)$ stated in the Appendix B.
It thereby follows that 
\begin{equation}
|N_-(\al_3)|^2=|M_{\al_2\al_1}^{\al_3;\ka_3}G_b(Q-2\al_3)|^2. 
\end{equation}
One has $|M_{\al_2\al_1}^{\al_3;\ka_3}|^2=e^{\pi i Q(Q-2\al_3)}$,
and $|G_b(Q-2\al_3)|^2=e^{-\pi i Q(Q-2\al_3)}|S_b(2\al_3)|^{-2}$ 
from the connection between $S_b$ and $G_b$, as well as the reflection
property of $S_b$ (see Appendix B). 
Therefore $|N_{-}(\al_3)|^2=|S_b(2\al_3)|^{-2}$.

The pole at $x=x_+$ corresponds to the pole at $y=0$ of $\Psi_b(U,V;W;y)$.
One may determine the singular term for $y\ra 0$ by applying 
Lemma \ref{doubpole} to the Euler integral representation \rf{Eulerint}
for the function $\Psi_b$: 
\begin{equation}
2\pi e^{-2\pi i y \be}\frac{G_b(-y+\ga-\be)}{G_b(\al)G_b(-y+Q)}
= \frac{1}{y} \frac{G_b(\ga-\be)}{G_b(\al)}+({\rm contributions\;\,regular\;\;
as}\;\;y\ra 0).
\end{equation}
The rest of the calculation proceeds as in the case of 
$N_-(\al_3)$ and yields $|N_{+}(\al_3)|^2=|S_b(2\al_3)|^{-2}$.
\end{proof}

\begin{propn}
Assume that the generalized eigenfunctions $\TP_{\al_3}$ 
decay exponentially for $\om\ra\infty$ and have
asymptotic behavior of the form 
\rf{TPasym} with $|N_{+}(\al_3)|^2=|N_{-}(\al_3)|^2$
for $\om\ra-\infty$. 
In that case one may define the ``inner product'' 
$(\Phi_{\al_3},\Phi_{\al_3'})$
as a bi-distribution which is explicitly given by
\begin{equation}
(\Phi_{\al_3},\Phi_{\al_3'})\;=\;|N_+(\al_3)|^2 \de(\al_3-\al_3').
\end{equation}
\end{propn}
\begin{proof}
Consider 
\begin{equation}\label{cal1}\begin{aligned}
( & C_{21}(\ka_3) \Phi_{\al_3},\Phi_{\al_3'})- 
(\Phi_{\al_3},C_{21}(\ka_3)\Phi_{\al_3'})
=\\
= & \lim_{W\ra\infty} \;\sum_{s=\pm}\;\int\limits_{-W}^{W}d\om\;\Bigl(
\bigl(\tilde{\de}_s(\om)\TP_{\al_3}(\om+sib)\bigl)^*\TP_{\al_3'}(\om)-
\bigl(\TP_{\al_3}(\om)\bigr)^*\tilde{\de}_s(\om)\TP_{\al_3'}(\om+sib)\Bigr),
\end{aligned}
\end{equation}
where the Fourier-transform of the explicit expression \rf{Cas} for 
$C_{21}(\ka_3)$ has been used. The contour of integration for the
second term in \rf{cal1} can be deformed into $\BR-isb$ plus
contours from $-W$ to $-W-isb$ and $W-isb$ to $W$. The integral 
over $\BR-isb$ cancels the first term on the right hand side of \rf{cal1}.
Only the contour from $-W$ to $-W-isb$ will give nonvanishing 
contributions in the limit $W\ra \infty$ due to the exponential decay 
of $\TP_{\al_3}(\om)$ for $\om\ra \infty$.
In the remaining term one gets in the limit $W\ra \infty$ contributions
only from the leading terms in the asymptotics of 
$\TP _{\al_3}(\om)$ for $\om\ra -\infty$ as quoted in Lemma \ref{TPasym}.
Taking into account that
\begin{equation}
\tilde{\de}_s(\om)=  \frac{1}{(q-q^{-1})^2}
e^{s\pi i b(Q-\al_1-\al_2)}+\CO(e^{2\pi b\om})
\end{equation}
for $\om\ra -\infty$, it follows that ($\al_3=\frac{Q}{2}+ip_3$, $\al_3'=\frac{Q}{2}+ip_3'$)
\begin{equation}\label{cal2}\begin{aligned}
(C_{21} (\ka_3) & \Phi_{\al_3},\Phi_{\al_3'})- 
(\Phi_{\al_3},C_{21}(\ka_3)\Phi_{\al_3'})
=\\
= & \frac{1}{(q-q^{-1})^2}
\lim_{W\ra\infty} \;\;\sum_{s=\pm}\;\sum_{\ep_1,\ep_2=\pm}
\frac{\bigl(N_{\ep_1}(\al_3)\bigr)^*N_{\ep_2}(\al_3')}
{2\pi i (\ep_1p_3-\ep_2p_3')}e^{2\pi i W(\ep_1p_3-\ep_2p_3')}\cdot\\
& \qquad\qquad\qquad\qquad\qquad\qquad\qquad\qquad\qquad
\cdot e^{2\pi s\ep_2 bp_3'}\bigl(1-e^{2\pi sb(\ep_1p_3-\ep_2p_3')}\bigr).
\end{aligned}
\end{equation}
The expression on the right hand side of \rf{cal2} vanishes by the
Riemann-Lebesque Lemma for $p_3\neq p_3'$ as well
as $\ep_1\neq \ep_2$.
The remainder is found to be
\begin{equation}\label{cal3}\begin{aligned}
(C_{21}(\ka_3) & \Phi_{\al_3},\Phi_{\al_3'})- 
(\Phi_{\al_3},C_{21}(\ka_3)\Phi_{\al_3'})
=\\
= & \;\bigl( [ip_3']^2_b-[ip_3]^2_b\bigr)\;|N_+(\al_3)|^2\;
\lim_{W\ra\infty} 
\frac{e^{2\pi i W(p_3-p_3')}-e^{-2\pi i W(p_3-p_3')}}{2\pi i (p_3-p_3')}.
\end{aligned}
\end{equation}
It follows that 
\begin{equation}\begin{aligned}
(\Phi_{\al_3},\Phi_{\al_3'})\;=& 
\;|N_+(\al_3)|^2
\lim_{W\ra\infty}\frac{e^{2\pi i W(p_3-p_3')}-
e^{-2\pi i W(p_3-p_3')}}{2\pi i (p_3-p_3')}\\ 
=& 
\;|N_+(\al_3)|^2 \;\de(\al_3-\al_3')
\end{aligned}\end{equation}
by the corresponding well-known property of the kernel $\sin(Rx)/x$,
cf. e.g. \cite[Chapter IX, Exercise 14]{RS2}.
\end{proof} 

\subsection{Intertwining property}

\begin{propn} \label{intertwprop}
The projections $\Pi_{21}(\al_3)$, $\al_3\in\BS$ map 
$\CP_{\al_2}\ot\CP_{\al_1}$ 
into $\CP_{\al_3}$ and 
satisfy the intertwining property \rf{intertw}.
\end{propn}
\begin{proof}
$F_f(\al_3,x_3)$ will be entire analytic w.r.t.
$x_3$ by straightforward application of Lemma \ref{doubpole},
using that $f$ is entire analytic in $x_2$, $x_1$ and the analytic
properties of the Clebsch-Gordan coefficients summarized in 
Lemma \ref{anasC}, Appendix C. One similarly finds by using
Lemma \ref{anasFTC}, Appendix C that the Fourier-transform
$F_f(\al_3,\ka_3)$ will be meromorphic in $\ka_3$ with poles
at $\ka=\pm(Q-\al+nb+mb^{-1})$, $n,m\in\BZ^{\geq 0}$
for any $\in\CP_{\al_2}\ot\CP_{\al_1}$. This establishes the
first claim in Proposition \ref{intertwprop}.

Note that the analytic continuation 
of the integral \rf{CGmap} that defines $F_f(\al_3,x_3)$
can be represented by 
integrating over a deformed contour $C^{(2)}\subset\BC^2$. 
For later use we will present suitable contours for the 
cases of 
analytic continuation to $\{x_3\in\BC;\Im(x_3)\in {[}0,\frac{b}{2}{]}\}$
and  $\{x_3\in\BC;\Im(x_3)\in[-\frac{b}{2},0]\}$ respectively:
In the first case one may integrate $x_1$ over the real axis
and instead of integrating over $x_2$ one may
integrate $x_{32}\equiv -iy_{32}$, cf. \rf{CGcoeff}, 
over a contour consisting of the union of
the half axes $(-\infty,-\de]$
and $[\de,+\infty)$, $b>\de>b/2$ 
with a half-circle  
in the upper half plane around $x_{32}=0$ 
of radius $\de$.
In the second case one may integrate $x_2$
over $\BR$, and $x_{31}\equiv -iy_{31}$ over
the contour $C_1$ consisting of the union of
the half axes $(-\infty,-\de]$
and $[\de,+\infty)$ with
a half-circle of radius $\de$ 
in the {\it lower} half plane around $x_{31}=0$.

Now consider the right hand side of \rf{intertw}.
The expressions for $\pi_{21}(u)$, $u=E,F,K$ contain the shift operators
\begin{equation}\label{shlist}
T_{x_1}^{+\frac{ib}{2}}T_{x_2}^{+\frac{ib}{2}},\quad
T_{x_1}^{-\frac{ib}{2}}T_{x_2}^{-\frac{ib}{2}}\quad\text{and}\quad
T_{x_1}^{-\frac{ib}{2}}T_{x_2}^{+\frac{ib}{2}}.
\end{equation}
The shift operator $T_{x_i}^{\pm\frac{ib}{2}}$ is ``partially integrated''
by (i) shifting the contour of integration over $x_i$ to the axis 
$\BR\mp\frac{ib}{2}$, where one will pick up a residue
contribution from the pole of the Clebsch-Gordan coefficients that lies between
these two contours, and (ii) introducing the new variables of integration
$x_i'\equiv x_i\pm\frac{ib}{2}$. In this way one rewrites the expression 
for $\CC_{21}\pi_{21}(u)f$ in the form
\begin{equation}\label{interst1}
\int\limits_{C_1}dx_2\int\limits_{C_2}dx_1\;
\Bigl(\pi_{21}^t(u) \CGC{\al_3}{x_3}{\al_2}{x_2}{\al_1}{x_1}\Bigr)\;
f(x_2,x_1),
\end{equation} 
where the $\pi_{21}^t$ denotes the transpose of $\pi_{21}$, and
the contours $C_i$, $i=1,2$ are just the contours introduced above 
to represent the analytic continuation w.r.t. $x_3$.
It is important to notice that
due to the fact that only the shift operators \rf{shlist}
appear in the expressions for 
$\pi_{21}(u)$, $u=E,F,K$ one does not need to introduce 
further deformations of the contours in order to treat
the poles from the factor in
the Clebsch-Gordan coefficients that depends on $x_2-x_1$ only.

It is verified by a straightforward calculation using \rf{S:funrel}
that the Clebsch-Gordan coefficients satisfy the finite difference 
equations
\begin{equation}
\pi_{21}^t(u) \CGC{\al_3}{x_3}{\al_2}{x_2}{\al_1}{x_1}=
\pi_{\al_3}(u)\CGC{\al_3}{x_3}{\al_2}{x_2}{\al_1}{x_1}, \qquad
u=E,F,K.
\end{equation}
Inserting these relations into \rf{interst1} yields an expression that is
easily identified as $\pi_{\al_3}(u)\CC_{21}f$. 
\end{proof}

\section{Racah-Wigner coefficients for $\USL$}

\subsection{Canonical decompositions for triple tensor products}

Triple tensor products $\CP_{\al_3}\ot
\CP_{\al_2}\ot\CP_{\al_1}$ carry a representation 
$\pi_{321}$ of $\USL$ given by 
\begin{equation} \begin{aligned}
\pi_{321}\;\equiv\;& (\pi_{\al_3}\ot
\pi_{\al_2}\ot\pi_{\al_1})\circ \De^{(3)},\\
\De^{(3)}\;\equiv\;&
(\De\ot\id)\circ\De\;\equiv\;(\id\ot\De)\circ\De.
\end{aligned}\end{equation}
The decomposition of this representation into
irreducibles can be constructed by iterating Clebsch-Gordan maps:
There are two canonical ways to do so, which will be referred to 
as ``s-channel'' and ``t-channel'' respectively. The first of these
corresponds to first decomposing the factor $\CP_{\al_2}\ot\CP_{\al_1}$
into a direct sum of irreducible representations $\CP_{\al_s}$
then performing the Clebsch-Gordan decomposition
of $\CP_{\al_3}\ot\CP_{\al_s}$. This extends to a unitary map 
\begin{equation}\label{CGmapis}
\CC_{3(21)}\quad : \quad  \begin{aligned}
{} & L^2(\BR\ti\BR\ti\BR)
\quad   \ra \quad 
L^2(\BS^2\ti\BR,d\mu(\al_4)d\mu(\al_s)dx_4)\\
& f(x_3,x_2,x_1) \quad \qquad  \ra \quad F_f^s(\al_4,\al_s,x_4), 
\end{aligned}
\end{equation}
The 
generalized Fourier-transform $F_f^s$ of $f$ is defined as 
\begin{equation}\label{21trsf}\begin{aligned}
F_f^s(\al_4,\al_s;x_4) \;\equiv\; \lim_{\ep_2\downarrow 0}
\lim_{\ep_1\downarrow 0}\int\limits_{\BR^2}
dx_3dx_{s}  & \;\,\CGC{\al_4}{x_4}{\al_3}{x_3}{\al_s}{x_{s}}^{}_{\ep_2}
\times \\ &\ti
\int\limits_{\BR^2}
dx_2dx_1\;\,\CGC{\al_s}{x_{s}}{\al_2}{x_2}{\al_1}{x_1}^{}_{\ep_1}\;
f(x_3,x_2,x_1),
\end{aligned}
\end{equation}
which in the notation $\fx\equiv(x_3,x_2,x_1)$, $d\fx\equiv
dx_3dx_2dx_1$ can be rewritten as 
\begin{equation}\label{Phis}\begin{aligned}
F_f^s(\al_4 & ,\al_s;x_4) \;\equiv\;\,  
\lim_{\ep\downarrow 0}\;\,\int\limits_{\BR^3}
d\fx\;\CpBls{\al_s}{\al_3}{\al_2}{\al_4}{\al_1}^{}_{\ep}(x_4;\fx)
\; f(\fx),\\
\text{where}\quad &  
\CpBls{\al_s}{\al_3}{\al_2}{\al_4}{\al_1}^{}_{\ep}(x_4;\fx)= 
\int_{\BR}dx_{s}\;\,
\CGC{\al_4}{x_4}{\al_3}{x_3}{\al_s}{x_{s}}^{}_{\ep}
\CGC{\al_s}{x_{s}}{\al_2}{x_2}{\al_1}{x_1}^{}_{\ep}
\quad\al_4,\al_s\in\BS,\;\,x_4\in\BR.
\end{aligned}\end{equation}
Some useful properties of the functions 
$\CpBls{\al_s}{\al_3}{\al_2}{\al_4}{\al_1}^{}_{\ep}(x_4;\fx)$ are 
collected in Appendix C.

The generalized Fourier-transformation $\CC_{3(21)}$ 
is such that the two-parameter family of projections
$\Pi^s(\al_4,\al_s):\CP_{\al_3}\ot\CP_{\al_2}\ot\CP_{\al_1}
\ra \CP_{\al_4}(\BR)$ defined 
by $f\ra F_f^s(\al_4,\al_s;.)$ 
intertwine the representation $\pi_{321}$
 with the irreducible representation $\pi_{\al_4}$.
It therefore realizes the following isomorphism of $\USL$ representations
\begin{equation}
\CP_{\al_3}\ot\CP_{\al_2}\ot\CP_{\al_1}\;\,\simeq\;\,
\int\limits_{\BS}^{\oplus}d\mu(\al_4)\;\,\CP_{\al_4}\,\ot\,\CS_{\mu}^{},
\end{equation}
where the multiplicity space $\CS_{\mu}\simeq
L^2(\BS,d\mu)$ is considered to be equipped with 
the trivial action of $\USL$.

A second canonical 
decomposition of $\CP_{\al_3}\ot\CP_{\al_2}\ot\CP_{\al_1}$ 
is obtained by first decomposing the factor $\CP_{\al_3}\ot\CP_{\al_2}$
into a direct sum of irreducible representations $\CP_{\al_t}$ and 
then performing the Clebsch-Gordan decomposition
of $\CP_{\al_t}\ot\CP_{\al_1}$. One obtains a map
\begin{equation}\label{CGmapit}
\CC_{(32)1}\quad : \quad  \begin{aligned}
{} & L^2(\BR\ti\BR\ti\BR)
\quad   \ra \quad 
L^2(\BS^2\ti\BR,d\mu(\al_4)d\mu(\al_t)dx_4) \\
& f(x_3,x_2,x_1) \quad \qquad  \ra \quad F_f^t(\al_4,\al_t,x_4), 
\end{aligned}
\end{equation}
where $F_f^t$ is defined by a generalized Fourier-transform of the same form
as \rf{21trsf} but with $\Phi_{21}^s$ replaced by 
\begin{equation} \label{Phit}
\CpBlt{\al_t}{\al_3}{\al_2}{\al_4}{\al_1}^{}_{\ep}(x_4;\fx)= 
\int_{\BR}dx_{t}\;\,
\CGC{\al_4}{x_4}{\al_t}{x_{t}}{\al_1}{x_1}^{}_{\ep}
\CGC{\al_t}{x_{t}}{\al_3}{x_3}{\al_2}{x_2}^{}_{\ep}.
\quad\al_4,\al_t\in\BS,\;\,x_4\in\BR.
\end{equation}
As in the case of the s-channel, one has 
a corresponding two-parameter family of projections
$\Pi^s(\al_4,\al_s):\CP_{\al_3}\ot\CP_{\al_2}\ot\CP_{\al_1}\ra \CP_{\al_4}$ 
that intertwine the representation $\pi_{321}$ 
with the irreducible representation $\pi_{\al_4}$.

\begin{rem}
The unitarity of the maps $\CC_{3(21)}$ and $\CC_{(32)1}$
ensures {\it existence} of self-adjoint extensions for the
operators $\pi_{3(21)}(u)$, $\pi_{(32)1}(u)$, $u=E,F,K,Q$:
Simply take the image of the self-adjoint extensions on 
$L^2(\BS^2\ti\BR)$ under $\CC_{3(21)}^{-1}$ or
$\CC_{(32)1}^{-1}$.

However, it is not a priori clear that such self-adjoint
extensions are {\it unique}. In particular, it could be that
the self-adjoint extensions that are defined in terms
of the maps $\CC_{3(21)}$ and $\CC_{(32)1}$ are inequivalent.
This disturbing possibility will be excluded shortly.
\end{rem}

\subsection{Relation between $\CC_{3(21)}$ and $\CC_{(32)1}$}

It will be convenient to also consider the Fourier-transforms
$\CpBl{\al_s}{\al_3}{\al_2}{\al_4}{\al_1}^{}_{\ep}(k_4;\fx)$, $\flat=s,t$
that are defined as 
\begin{equation}
\CpBl{\al_s}{\al_3}{\al_2}{\al_4}{\al_1}^{}_{\ep}(k_4;\fx)\;=\;
\int_{\BR}dx_4\;\,e^{2\pi ik_4x_4}\;
\CpBl{\al_s}{\al_3}{\al_2}{\al_4}{\al_1}^{}_{\ep}(x_4;\fx). \quad
\end{equation}

Unitarity of the maps $\CC_{3(21)}$ and $\CC_{(32)1}$ allows us to relate
the transforms $F_f^s$ and $F_f^t$ by a transformation of the form 
\begin{equation}\label{strel}
F_f^{s}(\al_4,\al_s,k_4)\;=\;
\int\limits_{\BS^2}d\al_4'd\al_t\int\limits_{\BR}dk_4\;\,
\CK\bigl[ \,\ew \begin{array}{lll} 
  {\scriptstyle \al_4} \ppe
  & {\scriptstyle \al_s} \ppe & {\scriptstyle k_4} \ew \\[-1.7mm] 
  {\scriptstyle \al_4'} \ppe 
  & {\scriptstyle \al_t} \ppe & {\scriptstyle k_4'} \ew\end{array}
  \bigr]\;
F_f^{t}(\al_4',\al_t^{},k_4').
\end{equation} 
The distribution $\CK$ appearing in \rf{strel}
can be represented as 
\begin{equation}\label{Krepr}\begin{aligned}
\CK & \bigl[ \,\ew \begin{array}{lll} 
  {\scriptstyle \al_4} \ppe
  & {\scriptstyle \al_s} \ppe & {\scriptstyle k_4} \ew \\[-1.7mm] 
  {\scriptstyle \al_4'} \ppe 
  & {\scriptstyle \al_t} \ppe & {\scriptstyle k_4'} \ew\end{array}
  \bigr]
=\\
& =\;\lim_{\rho\ra\infty}\lim_{\ep\downarrow 0}\;\int\limits_{-\infty}^{\infty}
dx_2\int\limits_{-\rho}^{\rho}dx_3dx_1\;\,
\bigl(\CpBlt{\al_t}{\al_3}{\al_2}{\al_4'}{\al_1}^{}_{\ep}(k_4';\fx)\bigr)^*
\CpBls{\al_s}{\al_3}{\al_2}{\al_4}{\al_1}^{}_{\ep}(k_4;\fx).
\end{aligned}\end{equation} 

We will first prove
\begin{propn} \label{Kform}
The distribution $\CK$ is of the form
\begin{equation}
\CK\bigl[ \,\ew \begin{array}{lll} 
  {\scriptstyle \al_4} \ppe
  & {\scriptstyle \al_s} \ppe & {\scriptstyle k_4} \ew \\[-1.7mm] 
  {\scriptstyle \al_4'} \ppe 
  & {\scriptstyle \al_t} \ppe & {\scriptstyle k_4'} \ew\end{array}
  \bigr]\;=\;\de(\al_4-\al_4')\de(k_4-k_4')\;
K\bigl[ \,\ew \begin{array}{ll} 
  {\scriptstyle \al_4} \ppe
  & {\scriptstyle \al_s} \ew \\[-1.7mm] 
  {\scriptstyle k_4} \ppe 
  & {\scriptstyle \al_t} \ew\end{array}
  \bigr].
\end{equation}
\end{propn}\begin{proof}
This will be a consequence of the following result: $\CK$ satisfies
\begin{equation}\label{adj}\begin{aligned}
\Bigl(\bigl[\al_4 -\fr{Q}{2}\bigr]_b^2-
      \bigl[\al_4'-\fr{Q}{2}\bigr]_b^2\Bigr)& 
 \;\,\CK\bigl[ \,\ew \begin{array}{lll}
   {\scriptstyle \al_4} \ppe
  & {\scriptstyle \al_s} \ppe & {\scriptstyle k_4} \ew \\[-1.7mm] 
  {\scriptstyle \al_4'} \ppe 
  & {\scriptstyle \al_t} \ppe & {\scriptstyle k_4'} \ew\end{array}
  \bigr]\;=\;0\\
(k_4-k_4')
 & \;\,\CK\bigl[ \,\ew \begin{array}{lll}
  {\scriptstyle \al_4} \ppe
  & {\scriptstyle \al_s} \ppe & {\scriptstyle k_4} \ew \\[-1.7mm] 
  {\scriptstyle \al_4'} \ppe 
  & {\scriptstyle \al_t} \ppe & {\scriptstyle k_4'} \ew\end{array}
  \bigr]\;=\;0.
\end{aligned}\end{equation}

To see that \rf{adj} implies the claim, consider 
the simplified case of a distribution $T\in\CS'(\BR)$
that satisfies $Tf=0$, where $f$ is a function that vanishes only at $x_0$
and such that $fg\in\CS(\BR)$
if $g\in\CS(\BR)$. This distribution has support only at $x_0$. 
By Theorem V.11 of \cite{RS} one has $T=\sum_{n=0}^Na_n(x_0)
\pa_x^n\de(x-x_0)$. 
It is then easy to see that $Tf=0$
implies $a_n=0$ for $n\neq 0$. 
The generalization to the case at hand is clear.

To verify \rf{adj} one may note that the functions 
$\CpBl{\al_t}{\al_3}{\al_2}{\al_4}{\al_1}^{}_{\ep}(k_4;\fx)$, $\flat=s,t$
satisfy eigenvalue equations for the 
operators $Q_{321}\equiv \pi_{321}(Q)$ and 
$K_{321}\equiv\pi_{321}(K)$
up to an error of order $\CO(\ep)$. It follows that
\begin{equation}\label{partint}\begin{aligned}
\Bigl( \bigl[ & \al_4 - \fr{Q}{2}\bigr]_b^2-  
      \bigl[\al_4'-\fr{Q}{2}\bigr]_b^2\Bigr)
\CK\bigl[ \,\ew \begin{array}{lll} 
  {\scriptstyle \al_4} \ppe
  & {\scriptstyle \al_s} \ppe & {\scriptstyle k_4} \ew \\[-1.7mm] 
  {\scriptstyle \al_4'} \ppe 
  & {\scriptstyle \al_t} \ppe & {\scriptstyle k_4'} \ew\end{array}
  \bigr]
=\\
& =\;\lim_{\ep_1,\ep_2\downarrow 0}\;\lim_{\rho\ra\infty}\;
\int\limits_{\BR}dx_2\int\limits_{-\rho}^{\rho}dx_3dx_1
\Bigl(
\bigl(\CpBlt{\al_t}{\al_3}{\al_2}{\al_4'}{\al_1}^{}_{\ep_1}(k_4';\fx)\bigr)^*
\;\,
Q_{321}\CpBls{\al_s}{\al_3}{\al_2}{\al_4}{\al_1}^{}_{\ep_2}(k_4;\fx)\\
 & \quad\qquad\qquad\qquad\qquad\qquad\qquad -
\bigl(Q_{321}\CpBlt{\al_t}{\al_3}{\al_2}{\al_4'}{\al_1}^{}_{\ep_1}
(k_4';\fx)\bigr)^*\;\,
\CpBls{\al_s}{\al_3}{\al_2}{\al_4}{\al_1}^{}_{\ep_2}(k_4;\fx)\Bigr).
\end{aligned}\end{equation}
The right hand side of \rf{partint} will vanish if 
$Q_{321}$ can be ``partially integrated''. 
To show that this is the case, one needs
some information on the form that $Q_{321}$ takes when acting
on functions $f(\fx)$. By straightforward evaluation of its
definition one obtains an expression in terms of shift operators
\[
T_{1}^{is_1b}T_{2}^{is_2b}T_{3}^{is_3b},\quad\text{where}\;\; T_i=T_{x_i}, 
\;\; s_i\in\{+,-\},\;\;i=1,2,3.
\]
It is convenient to introduce an alternative set of shift operators
\[ 
T_+^3=T_1^{}T_2^{}T_3^{},\qquad T_{21}^2=T_2^{}T_1^{-1}\qquad
T_{32}^{2}=T_3^{}T_2^{-1}.
\]
The crucial point now is that the expression for $Q_{321}$ when 
rewritten in terms of $T_+$, $T_{21}$, $T_{32}$ takes the following
form
\begin{equation}
Q_{321}^{}\;=\;\sum_{n_+=-3}^3\sum_{n_{21}=0}^{3}\sum_{n_{32}=0}^{3}
P_{n_+n_{21}n_{32}}(\fx)\;T_+^{in_+b}\;T_{21}^{\frac{2}{3}ibn_{21}}
\;T_{32}^{\frac{2}{3}ibn_{32}},
\end{equation}
so it contains shifts of $x_{21}$, $x_{32}$, $x_{31}$ by {\it positive} 
imaginary amounts up to $2ib$ only. Furthermore note that
in \rf{partint} one may replace $T_+$ by $e^{-2\pi i k_4}$. 
The analytic properties of the integrand in \rf{partint} as following from
Lemma \ref{anas} in Appendix C now allow to partially integrate 
$Q_{321}$ by appropriate shifts of the contours of integration over 
$x_3,x_2,x_1$ (cf. proof of Proposition \ref{intertwprop}).

The verification of the second equation in \rf{adj} is similar.
\end{proof}

\begin{rem}
This result implies that the self-adjoint extensions of
$\pi_{321}(u)$, $u=K,Q$ 
that are defined by the maps 
$\CC_{3(21)}$ and $\CC_{(32)1}$ indeed coincide.
A similar argument as in the proof of the previous proposition
will also cover the two other cases $u=E,F$.
\end{rem}

\subsection{Calculation of the Racah-Wigner coefficients I} 

It will be useful to also introduce  \begin{equation}\label{Xrepr}
\begin{aligned}
\CX & \bigl[ \,\ew \begin{array}{lll} 
  {\scriptstyle \al_4} \ppe
  & {\scriptstyle \al_s} \ppe & {\scriptstyle x_4} \ew \\[-1.7mm] 
  {\scriptstyle \al_4'} \ppe 
  & {\scriptstyle \al_t} \ppe & {\scriptstyle x_4'} \ew\end{array}
  \bigr]
=\\
& =\;\lim_{\ep\ra 0+}\;\int\limits_{-\infty}^{\infty}
dx_3dx_2dx_1\;\,
\bigl(\CpBlt{\al_t}{\al_3}{\al_2}{\al_4}{\al_1}^{}_{\ep}(x_4';\fx)\bigr)^*\;
\CpBls{\al_s}{\al_3}{\al_2}{\al_4}{\al_1}^{}_{\ep}(x_4;\fx).
\end{aligned}\end{equation}
Proposition \ref{Kform} has an obvious counterpart for $\CX$:
\begin{propn} 
The distribution $\CX$ is of the form
\begin{equation}\label{Xform}
\CX\bigl[ \,\ew \begin{array}{lll} 
  {\scriptstyle \al_4} \ppe
  & {\scriptstyle \al_s} \ppe & {\scriptstyle x_4} \ew \\[-1.7mm] 
  {\scriptstyle \al_4'} \ppe 
  & {\scriptstyle \al_t} \ppe & {\scriptstyle x_4'} \ew\end{array}
  \bigr]\;=\;\de(\al_4-\al_4')\de(x_4-x_4')\;
\SJS{\al_1}{\al_2}{\al_3}{\al_4}{\al_s}{\al_t}.
\end{equation}
\end{propn}
\begin{proof}
Introduce 
\begin{equation}\label{regK}
\CK_{\ep,\rho}^{}  \bigl[ \,\ew \begin{array}{lll} 
  {\scriptstyle \al_4} \ppe
  & {\scriptstyle \al_s} \ppe & {\scriptstyle k_4} \ew \\[-1.7mm] 
  {\scriptstyle \al_4'} \ppe 
  & {\scriptstyle \al_t} \ppe & {\scriptstyle k_4'} \ew\end{array}
  \bigr]
=
\;\int\limits_{-\infty}^{\infty}dx_2\int\limits_{-\rho}^{\rho}dx_3dx_1\;\,
\bigl(\CpBls{\al_s}{\al_3}{\al_2}{\al_4}{\al_1}^{}_{\ep}(k_4;\fx)\bigr)^*
\CpBlt{\al_t}{\al_3}{\al_2}{\al_4}{\al_1}^{}_{\ep}(k_4;\fx).
\end{equation}
The coefficient of $\de(k_4-k_4')$ in the expression 
for $\CK$ coinicides with the sum of the coefficients with which
$e^{-2\pi i(k_4-k_4')x_1}$ and $e^{-2\pi i(k_4-k_4')x_3}$ appear
in the asymptotic expansion of the integrand in \rf{regK}, 
cf. Lemma \ref{anas}. Lemma \ref{distPW} identifies the origin of these
terms in the asymptotic expansion of $\Phi^{\flat}$, $\flat=s,t$, with 
the poles in the dependence of $\Phi^{\flat}[\ldots]_{\ep}(x_4;\fx)$, $\flat=s,t$ 
on their variable $x_4$. It follows that the coefficient of $\de(k_4-k_4')$
in the expression 
for $\CK$ is independent of $k_4$. The result now 
follows from standard
properties of the Fourier transformation.
\end{proof}
\begin{propn} \label{RWformprop} We have
\begin{equation}\label{RWformula}\begin{aligned}
\SJS{\al_1}{\al_2}{\al_3}{\al_4}{\al_s}{\al_t}\;\, = & \;\,N\;\,  
 \frac{S_b(\al_2+\al_s-\al_1)
S_b(\al_t+\al_1-\al_4)}
     {S_b(\al_2+\al_t-\al_3)
S_b(\al_s+\al_3-\al_4)}\cdot\\
&\cdot|S_b(2\al_t)|^2
\int\limits_{-i\infty}^{i\infty}ds \;\;
\frac{S_b(U_1+s)S_b(U_2+s)S_b(U_3+s)S_b(U_4+s)}
{S_b(V_1+s)S_b(V_2+s)S_b(V_3+s)S_b(V_4+s)},
\end{aligned} 
\end{equation}
where the coefficients $U_i$ and $V_i$, $i=1,\ldots,4$ are given by
\begin{equation}  
\begin{aligned} U_1=& \al_s+\al_1-\al_2 \\
        U_2=& Q+\al_s-\al_2-\al_1 \\
        U_3=& \al_s+\al_3-\al_4 \\
        U_4=& Q+\al_s-\al_3-\al_4
\end{aligned}\qquad
\begin{aligned} 
        V_1=& 2Q+\al_s-\al_t-\al_2-\al_4 \\
        V_2=& Q+\al_s+\al_t-\al_4-\al_2 \\
        V_3=& 2\al_s \\
        V_4=& Q,
\end{aligned}
\end{equation}
and $N$ is a constant.
\end{propn}
\begin{proof} Let
\begin{equation}\label{regX}
\CK_{\ep}^{}  \bigl[ \,\ew \begin{array}{lll} 
  {\scriptstyle \al_4} \ppe
  & {\scriptstyle \al_s} \ppe & {\scriptstyle x_4} \ew \\[-1.7mm] 
  {\scriptstyle \al_4'} \ppe 
  & {\scriptstyle \al_t} \ppe & {\scriptstyle x_4'} \ew\end{array}
  \bigr]
=
\;\int\limits_{-\infty}^{\infty}dx_3dx_2 dx_1\;\,
\bigl(\CpBlt{\al_t}{\al_3}{\al_2}{\al_4'}{\al_1}^{}_{\ep}(x_4';\fx)\bigr)^*
\;\CpBls{\al_s}{\al_3}{\al_2}{\al_4}{\al_1}^{}_{\ep}(x_4;\fx).
\end{equation}
The analytic and
asymptotic properties of the
integrand follow from Lemma \ref{anasX} in Appendix C. Let us observe that
for $\ep>0$ one is dealing with absolutely convergent integrals, 
the integrand being meromorphic both w.r.t. the integration variables
and the parameters. 
The integral \rf{regX} therefore does not depend on the order
in which the integrations are performed, so we will assume that
it is first integrated over $x_2$. 

Singular behavior will emerge in the limit $\ep\ra 0$.
We will call 
a pole relevant if it has distance of $\CO(\ep)$
from the real axis, irrelevant otherwise\footnote{We of course assume
that $\ep$ has been chosen to be much smaller than $b$}.
It then easily follows from
Lemma \ref{doubpole} that the integration 
over $x_2$ does not introduce any new relevant poles 
since all the relevant poles in the
$x_2$ dependence that have distance of $\CO(\ep)$ are lying on the 
same side of the contour.

Next one may integrate over $x_1$. We find from Lemma \ref{anasX} in 
Appendix C that
\begin{equation}\begin{aligned}
\CpBls{\al_s}{\al_3}{\al_2}{\al_4}{\al_1}^{}_{\ep}(x_4,\fx)\;=\;&
\frac{\CR_{13}^s}{x_1-x_3+\al_{13}-2i\ep}+
\frac{\CR_{14}^s}{x_1-x_4+\al_{14}-2i\ep}+(\text{Reg}_s),\\
\bigl(\CpBlt{\al_t}{\al_3}{\al_2}{\al_4'}{\al_1}^{}_{\ep}(x_4',\fx)\bigr)^*
\;=\;&
\frac{\CR_{13}^t}{x_1-x_3+\al_{13}'+2i\ep}+
\frac{\CR_{14}^t}{x_1-x_4'+\al_{14}'+i\ep}+(\text{Reg}_t),
\end{aligned}\end{equation}
where $(\text{Reg}_{\flat})$, $\flat=s,t$ are terms that 
do not lead to relevant poles in the variable $x_1$
after having integrated over $x_2$.
The following abbreviations have been used:
\begin{equation}\begin{aligned}
\al_{13}=& \fr{i}{2}(\al_1+\al_3-2(Q-\al_4)),\\
\al_{14}=& \fr{i}{2}(\al_1-\al_4),
\end{aligned}
\qquad
\begin{aligned}
\al_{13}'=& \fr{i}{2}(\al_1+\al_3-2(Q-\al_4')),\\
\al_{14}'=& \fr{i}{2}(\al_1-\al_4').
\end{aligned}
\end{equation}
It is then easily found by using 
Lemma \ref{doubpole} that the result of the integration over $x_1$
will have poles at the following locations:
\begin{equation}\label{fourpoles}\begin{aligned}
{} & i(\al_4-\al_4')-4i\ep=0,\\
& x_4'-x_4+\fr{i}{2}(\al_4'-\al_4)-3i\ep=0,
\end{aligned}\qquad
\begin{aligned}
{} & x_3-x_4-\fr{i}{2}(\al_3+\al_4-2(Q-\al_4'))-4i\ep=0,\\
& x_4'-x_3+\fr{i}{2}(\al_3+\al_4'-2(Q-\al_4))-3i\ep=0.
\end{aligned}
\end{equation}
The relevant residues can easily be assembled from the 
expressions given in Appendix C. Moreover, it is straightforward to
work out their poles. By again using
Lemma \ref{doubpole} one then finds
that all four poles listed in 
\rf{fourpoles} will, after doing the $x_3$ integration,
produce terms that are singular for $x_4=x_4'$, $\al_4=\al_4'$ and 
$\ep\ra 0$. The terms that lead to $\de(x_4-x_4')\de(\al_4-\al_4')$
are easily identified by means of
\begin{equation}
\lim_{\ep\ra 0+}\Bigl(\frac{1}{x-i\ep}-\frac{1}{x+i\ep}\Bigr)
=2\pi i \de(x).
\end{equation}
All these terms have as residue an expression proportional to
\begin{equation}\label{resform}\begin{aligned}
\Res_{y_{31}=0}\Res_{y_{21}=0}& 
\CGC{\al_4}{\ast}{\al_3}{\ast}{\al_s}{\ast}
 \Res_{y_{31}=0}\Res_{y_{21}=0}
\CGC{\al_4}{\ast}{\al_t}{\ast}{\al_1}{\ast}\\
& \int_{\BR}dx_2\;
\Res_{y_{31}=0}
\CGC{\al_s}{\ast}{\al_2}{x_2}{\al_1}{x_1}_{x_1=
x_3-\al_{13}}
\Res_{y_{32}=0}
\CGC{\al_t}{x_t}{\al_3}{\ast}{\al_2}{x_2}_{x_t=
x_3-\frac{i}{2}(\al_3-\al_t)}.
\end{aligned}\end{equation}
One just needs to assemble the ingredients to check that
the expression \rf{resform} coincides with what one finds on the
right hand side of \rf{RWformula} 
\end{proof}
\begin{rem}
With more patience, one could of course also fix the 
constant $N$ by the method used in the previous proof. 
We refrain from doing so since we will present a less tedious 
and more illuminating way of calculating it in the next subsection.
What will be needed there, however, is the information
on analyticity of the coefficients $\{\ldots\}$ w.r.t. $\al_t$ that
follows from Proposition \ref{RWformprop}. 
\end{rem}

\subsection{Relation between the distributions $\Phi^s$ and $\Phi^t$}

\begin{propn}\label{Phi_s-t_rel}
$\Phi^s$ and $\Phi^t$ are related by a linear 
transformation of the
form
\begin{equation}\label{Racahdef}
\CpBls{\al_{s}}{\al_3}{\al_2}{\al_4}{\al_1}(x_4;\fx)\;=\;
\int\limits_{\BS}d\al_{t}\;\,
\SJS{\al_1}{\al_2}{\al_3}{\al_4}{\al_{s}}{\al_{t}}\;\,
\CpBlt{\al_{t}}{\al_3}{\al_2}{\al_4}{\al_1}(x_4;\fx).
\end{equation}
The relation \rf{Racahdef} can be read either as (i) relation between 
function analytic in 
\[ \begin{aligned}
\CA^{(4)}\equiv 
\{ \ux=(x_4,x_3,x_2,x_1)\in\BC^4;\;\, & \Im(x_1)<\Im(x_2)<\Im(x_3),\\
& \Im(x_1)<\Im(x_4)<\Im(x_3), \;\Im(x_3-x_1)<Q\},
\end{aligned} \]
or (ii) as relation between functions meromorphic w.r.t. $\ux\in\BC^4$,
or (iii) as relation between distributions defined as boundary values 
of $\Phi^{\flat}$, $\flat=s,t$ for $(x_4,\fx)\in\BR^4$.
\end{propn}
\begin{proof}
We will start from equation \rf{strel}. 
By using Fourier-transformation w.r.t. the variable $k_4$ and  
equation \rf{Xform} one may rewrite \rf{strel} as follows:
\begin{equation}\label{strel2}
F_f^{s}(\al_4,\al_s,x_4)\;=\;
\int\limits_{\BS}d\al_t\;\,
\SJS{\al_1}{\al_2}{\al_3}{\al_4}{\al_s}{\al_t}\;
F_f^{t}(\al_4',\al_t^{},x_4).
\end{equation}
Let us introduce sequences of test-functions that tend towards
delta-distributions:
\begin{equation}
t_n(\fy;\fx)\;=\; \Bigl(\frac{n}{2\pi}\Bigr)^{\frac{3}{2}}e^{-\frac{n}{2}
|\!|\fx-\fy|\!|^2}, \qquad \fy=(y_3,y_2,y_1).
\end{equation}
\begin{lem}\label{dellim}
Let $\uy\equiv(x_4,\fy)\in\CA^{(4)}$ with $\Im(y_1)<0$. 
In this case one has 
\begin{equation}
\lim_{n\ra\infty}F^{\flat}_{t_n(\fy;.)}(\al_4,\al_{\flat},x_4)
=\CpBl{\al_{\flat}}{\al_3}{\al_2}{\al_4}{\al_1}(x_4;\fy).
\end{equation}
\end{lem}
\begin{proof}
By writing out the definition of $F^{\flat}_{t_n}$ and shifting the contours
of integration over $x_i$ to $\BR+i\Im(y_i)$, $i=1,2,3$, 
one reduces the claim to the standard result that
\[
\lim_{n\ra\infty}t_n(\fy;\fx)=\de^3(\fx-\fy)
\]
for $\Im(y_i)=0$, $i=1,2,3$ (Note that $\Phi^{\flat}$ is regular
for these values  of its arguments as follows from
Lemma \ref{anasX}, Appendix C). \end{proof}
We will now consider the sequence with elements
\begin{equation}\label{seq}
\int\limits_{\BS}d\al_t\;\,
\SJS{\al_1}{\al_2}{\al_3}{\al_4}{\al_s}{\al_t}\;
F_{t_n(\fy,.)}^{t}(\al_4,\al_t,x_4).
\end{equation}
It converges for $n\ra\infty$
due to Lemma \ref{dellim} and equation \rf{strel2}. 
We would like to show that one may
exchange the limit $n\ra\infty$ with the integration over $\al_t$ so that 
the limit of \rf{seq} is given by the integral
\begin{equation}\label{seqlim}
\int\limits_{\BS}d\al_t\;\,
\SJS{\al_1}{\al_2}{\al_3}{\al_4}{\al_s}{\al_t}\;
\CpBlt{\al_{t}}{\al_3}{\al_2}{\al_4}{\al_1}(x_4;\fy).
\end{equation}
To this aim it is useful to note 
that 
\begin{lem} \label{expdec}
Under the conditions on the variable $\fy$ introduced in Lemma \ref{dellim}
one finds that the integrand in \rf{seqlim}
decays exponentially for
$p_t\equiv -i(\al_t-\frac{Q}{2})\ra \pm \infty$. The integrand in 
\rf{seq} decays at least as fast as the integrand in \rf{seqlim}.
\end{lem}
\begin{proof}
By a straightforward 
calculation using the method in the proof of Lemma \ref{asymeuler}, 
Appendix B and eqn. \rf{SGas} one finds that 
\begin{equation}\begin{aligned}
{} & 
\CpBlt{\al_{t}}{\al_3}{\al_2}{\al_4}{\al_1}(x_4;\fy)\quad\text{decays stronger
than $e^{\mp\pi Qp_t}$ and}\\
 & 
\SJS{\al_1}{\al_2}{\al_3}{\al_4}{\al_s}{\al_t}
\quad\text{grows as
$e^{\pm\pi Qp_t}$}
\end{aligned}\end{equation}
for $p_t\ra\infty$. The first statement in Lemma \ref{expdec} follows.

The second statement follows from 
the first by shifting the contour of integration over $x_1$ in the 
definition of $F^t_{t_n(\fy,.)}$ to $\BR+i\Im(y_1)$.
\end{proof}
The integrals \rf{seq}\rf{seqlim} can therefore be transformed into
integrals over a compact set, e.g. the interval $[0,1]$. In order to justify 
the exchange of limit and integration it therefore suffices to prove 
the following 
\begin{lem}\label{unif}
The convergence of $F_{t_n(\fy,.)}^{t}(\al_4,\al_t,x_4)$ is 
uniform in $\al_t$.
\end{lem}
\begin{proof} To shorten the exposition, let us consider a slightly 
simplified situation. Assume that $f_{p}(x)$ is analytic w.r.t. both
$p$ and $x$ in open strips that contain the real axis and decays
exponentially for either $|p|$ or $|x|$ going to infinity. Let 
$t_n(x)=\sqrt{\frac{n}{2\pi}}e^{-nx^2/2}$ and study the convergence
of $f_{p,n}\equiv \int_{\BR} dx f_{p}(x)t_n(x)$ for $n\ra \infty$.
Upon writing $f_p(x)=f_p(0)+xg_p(x)$, the task reduces to the study
of 
\begin{equation}\label{deerr}
  \int\limits_{\BR} dx \;\, g_{p}(x) \;\, x t_n(x)\;=\;\frac{1}{\sqrt{2\pi n}}
\int\limits_{\BR} dx\;\,e^{-\frac{n}{2}x^2}\;
\pa_xg_p(x).  
\end{equation}
Convergence for $n\ra\infty$ will be uniform in $p$ provided that 
$\pa_x g_p(x)$ is bounded as function of both $p$ and $x$. But 
this is a consequence of our assumptions: The exponential decay
allows us to transform $f_p(x)$ (resp. $\pa_xg_p(x)$) to a function  
that is analytic on a compact rectangle in $\BC^2$, and therefore bounded.

The regularity properties of $\Phi^t$  necessary to extend the 
argument to the present situation follow 
from Lemma \ref{anasX}, Appendix C.\end{proof}
We have proved \rf{Racahdef} provided $(x_4,\fx)$ satisfies the same conditions
as $(x_4,\fy)$ in Lemma \ref{dellim}. Proposition \ref{Phi_s-t_rel}
follows by analytic
continuation.  
\end{proof}

\subsection{Calculation of Racah-Wigner coefficients II} \label{RWII}

We have shown that the meromorphic functions $\Phi^s$ and $\Phi^t$
are related by an integral transformation of the form
\rf{Racahdef}.
If one fixes the values of 
three of the four variables $x_4,\ldots,x_1$ in \rf{Racahdef} one obtains
an integral transformation for a function of a single variable. In fact, the
analytic properties of $\Phi_{\al_s}^s$ and $\Phi_{\al_t}^t$ even allow
one to choose complex values. It will be convenient to consider
\begin{equation} \label{limits}
\CBls{\al_s}{\al_3}{\al_2}{\bar{\al}_4}{\al_1}{\Psi}(x)=
\lim_{x_4\ra\infty}e^{2\pi\al_4x_4}
\lim_{x_2\ra-\infty}
\prod_{j=1}^{3}e^{-2\pi \al_jx_j}
\CBls{\al_s}{\al_3}{\al_2}{\bar{\al}_4}{\al_1}{\Phi}(\fx)
\Bigr|_{x_3=\frac{i}{2}(Q+\al_2-\al_4)}^{x_1=x},
\end{equation} 
where $\bar{\al}=Q-\al$, 
and the same for $\Psi^t_{\al_t}$. The integral that defines 
$\Phi_{\al_s}^s$ and $\Phi_{\al_t}^t$, \rf{Phis}\rf{Phit} 
can be done explicitly in this limit by using \rf{Eulerint}. 
One finds expressions of the form
\begin{equation}\begin{aligned}
{} & \CBls{\al_s}{\al_3}{\al_2}{\bar{\al}_4}{\al_1}{\Psi}(x)= 
\CBls{\al_s}{\al_3}{\al_2}{\bar{\al}_4}{\al_1}{N}
\CBls{\al_s}{\al_3}{\al_2}{\bar{\al}_4}{\al_1}{\Theta}(x)\\
& \qquad \CBls{\al_s}{\al_3}{\al_2}{\bar{\al}_4}{\al_1}{\Theta}(x)=
e^{+2\pi x(\al_s-\al_2-\al_1)}
F_b(\al_s+\al_1-\al_2,\al_s+\al_3-\al_4;2\al_s;-ix)\\
 & \CBlt{\al_t}{\al_3}{\al_2}{\bar{\al}_4}{\al_1}{\Psi}(x)= 
\CBlt{\al_t}{\al_3}{\al_2}{\bar{\al}_4}{\al_1}{N}
\CBlt{\al_t}{\al_3}{\al_2}{\bar{\al}_4}{\al_1}{\Theta}(x)\\
& \qquad\CBlt{\al_t}{\al_3}{\al_2}{\bar{\al}_4}{\al_1}{\Theta}(x)=
e^{-2\pi x(\al_t+\al_1-\al_4)}
F_b(\al_t+\al_3-\al_2,\al_t+\al_1-\al_4;2\al_t;+ix),
\end{aligned}
\end{equation}
where $F_b$ is the b-hypergeometric function defined in the Appendix,
and $N_{\al_s}^s$, $N_{\al_t}^t$ are certain normalization factors.

The linear transformation following
from \rf{Racahdef} can now be calculated as follows: One observes that
$\Psi^s_{\al_s}$ (resp. $\Psi^t_{\al_t}$) are eigenfunctions
of the finite difference operators $\CQ_{s}$ and $\CQ_{t}$
defined respectively by
\begin{equation}\begin{aligned}
\CQ_s= & \bigl[d_x+\al_1+\al_2-\fr{Q}{2}\bigr]^2-
e^{+2\pi b x}\bigl[d_x+\al_1+\al_2+\al_3-\al_4\bigr]
\bigl[d_x+2\al_1\bigr]\\
 \CQ_t= & \bigl[d_x+\al_1-\al_4+\fr{Q}{2}\bigr]^2-
e^{-2\pi b x}\bigl[d_x+\al_1+\al_2-\al_3-\al_4\bigr]\bigl[d_x\bigr].
\end{aligned}\end{equation}
It can be shown that 
\begin{thm} The operators $\CQ_s$ and $\CQ_t$ have unique self-adjoint 
extensions in $L^2(\BR, dxe^{2\pi Qx})$. Bases of $L^2(\BR, dxe^{2\pi Qx})$ 
in the sense of
generalized eigenfunctions are given by the sets of functions
$\{ \Theta^s_{\al_s};\al_s\in\BS\}$ and 
$\{ \Theta^t_{\al_t};\al_t\in\BS\}$, where the normalization is given by 
\begin{equation}
\int\limits_{\BR} dx \;e^{2\pi Qx}\; \bigl(
\CBlf{\al_{\flat}'}{\al_3}{\al_2}{\bar{\al}_4}{\al_1}{\Theta}(x) \bigr)^* 
\CBlf{\al_{\flat}}{\al_3}{\al_2}{\bar{\al}_4}{\al_1}{\Theta}(x)=
\de(\al_{\flat}-\al_{\flat}'), \quad\flat=s,t.
\end{equation}\end{thm}
The proof is omitted as it is very similar to the proof
of Theorem \ref{casspec}.
It follows that the Racah-Wigner coefficients can be evaluated
in terms of the overlap between these two bases:
\begin{equation}
\SJS{\al_1}{\al_2}{\al_3}{\bal_4}{\al_s}{\al_t}\;=\;
\frac{\CBls{\al_s}{\al_3}{\al_2}{\bar{\al}_4}{\al_1}{N}}
{\CBlt{\al_t}{\al_3}{\al_2}{\bar{\al}_4}{\al_1}{N}}
\;\,\int\limits_{\BR} dx \;e^{2\pi Qx}\; \bigl(
\CBlt{\al_t}{\al_3}{\al_2}{\bar{\al}_4}{\al_1}{\Theta}(x) \bigr)^* 
\CBls{\al_s}{\al_3}{\al_2}{\bar{\al}_4}{\al_1}{\Theta}(x).
\end{equation}
The integral can be done by using the representation \rf{Barnesint}
for the b-hypergeometric function. The result is
just equation \rf{RWformula} with $N=1$.

\subsection{Properties the Racah-Wigner coefficients}

First of all let us note that orthogonality and completeness of the bases 
$\{ \Phi^s_{\al_s};\al_s\in\BS\}$ and 
$\{ \Phi^t_{\al_t};\al_t\in\BS\}$ imply
the following orthogonality relations for the b-Racah-Wigner symbols
\begin{equation}\label{orth3}
\int\limits_{\BS}d\al_s\;|S_b(2\al_s)|^2 
\; \SJS{\al_1}{\al_2}{\al_3}{\al_4}{\al_s}{\al_t}
\bigl(\SJS{\al_1}{\al_2}{\al_3}{\al_4}{\al_s}{\be_{t}}\bigr)^*\;=\;
|S_b(2\al_t)|^2\;\de(\al_t-\be_{t}).
\end{equation}
This may be verified e.g. by 
rewriting 
\begin{equation}\begin{aligned}\label{t-norm}
\bigl(\,\CpBlt{\al_t}{\al_3}{\al_2}{\al_4}{\al_1}^{}_{\ep}(x_4;.) & \,,\,
\CpBlt{\al_t'}{\al_3}{\al_2}{\al_4'}{\al_1}^{}_{\ep}(x_4';.)\,\bigr)\;=\\
& =\;
|S_b(2\al_t)|^{-2}\de(\al_t-\al_t')\de(\al_4-\al_4')\de(x_4-x_4')
\end{aligned}\end{equation}
with the help of the inversion formula to \rf{Racahdef}
\begin{equation}
\CpBlt{\al_{s}}{\al_3}{\al_2}{\al_4}{\al_1}(x_4;\fx)\;=\;
\int\limits_{\BS}d\al_{s}\;\biggl|\frac{S_b(2\al_s)}{S_b(2\al_t)}\biggr|^2\;
\bigl(\SJS{\al_1}{\al_2}{\al_3}{\al_4}{\al_{s}}{\al_{t}}\bigr)^*
\;\CpBls{\al_{t}}{\al_3}{\al_2}{\al_4}{\al_1}(x_4;\fx),
\end{equation}
and finally using \rf{t-norm} with subscripts $t$ replaced by $s$.

Second, by considering quadruple products of representations
one finds the so-called pentagon equation in the usual way:
\begin{equation}
\int\limits_{\BS}d\de_1\;
\SJS{\al_1}{\al_2}{\al_3}{\be_2}{\be_1}{\de_1}
\SJS{\al_1}{\de_1}{\al_4}{\al_4}{\be_2}{\ga_2}
\SJS{\al_2}{\al_3}{\al_4}{\ga_2}{\de_1}{\ga_1}
\;\,=\;\, 
\SJS{\be_1}{\al_3}{\al_4}{\al_5}{\be_2}{\ga_1}
\SJS{\al_1}{\al_2}{\ga_1}{\al_5}{\be_1}{\ga_2}.
\end{equation}

\subsection{From intertwiners to coinvariants}

Let us consider coinvariants on tensor products of representations.
These will be maps $\CB:
\CP_{\al_n}\ot\ldots\ot\CP_{\al_1.}
\ra \BC$ that satisfy the coinvariance property 
\begin{equation}\label{invprop}
\CB\circ\bigl((\pi_{\al_n}\ot\ldots\ot\pi_{\al_1})\De^{(n)}(u)\bigl)\;=\;0,
\quad u\in\USL,
\end{equation}
where $\De^{(n)}$ is defined recursively by $\De^{(n)}=(\id\ot\De)
(\De^{(n-1)})=
(\De\ot\id)(\De^{(n-1)})$, $\De^{(2)}\equiv\De$.

The basic case to consider is $n=2$. Let
$\CB_{\al}:\CP_{Q-\al}\ot\CP_{\al} 
\ra \BC$ be defined by 
\begin{equation}
\CB_{\al}(f\ot g)\;\equiv\; \bra \,f\,, \, 
\CT g\,\ket,\qquad \CT\;\equiv\; T_x^{-i\frac{Q}{2}}
\end{equation}
\begin{propn}
$\CB_{\al}$ satisfies the coinvariance property \rf{invprop}.
\end{propn}
\begin{proof}
Let us note that 
\begin{equation}
\bra \, T_x^{i\al}f\,,\,g\,\ket\;=\; \bra \, f\,,T_x^{-i\al}g\,\ket
\end{equation}
if $f\in\CP_{Q-\al}$ and $g\in\CP_{\al}$. 
A straightforward calculation then shows that
\begin{equation}\label{transp}
\bra \, \pi_{Q-\al}(u)f\,,\,g\,\ket\;=\; \bra \, f\,,\pi_{\al}(u)g\,\ket,
\quad u\in\USL.
\end{equation}
It is useful to also note the commutation relations 
\begin{equation}
\CT\,E_{\al}\;=\; e^{-i\pi bQ}\,E_{\al}\,\CT,\qquad
\CT\,F_{\al}\;=\; e^{+i\pi bQ}\,F_{\al}\,\CT,\qquad
\CT\,K_{\al}\;=\; K_{\al}\,\CT.
\end{equation}
We may then calculate in the case $u=E$
\begin{equation}\begin{aligned}
\CB_{\al}\bigl( ((\pi_{Q-\al} & \ot\pi_{\al})\circ\De(E)) f\ot g\bigr)\; 
 =\\
&=\;
\bra \, E_{Q-\al}f\,,\,\CT K_{\al}g\,\ket+
\bra \, K_{Q-\al}f\,,\,\CT E_{\al}g\,\ket \\
&=\; \bra \, E_{Q-\al} f\,, \, K_{\al}\CT g\,\ket+
e^{-i\pi bQ} \bra \,K_{Q-\al} f\,,\,E_{\al}\CT g\,\ket\\
&=\; \bra \, f\,,\, E_{\al}K_{\al}\CT g\,\ket-
q^{-1} \bra \,\CT f\,,\,K_{\al}E_{\al}\CT g\,\ket\\
&=\;0.
\end{aligned}
\end{equation}
The calculation for the case $u=F$ is identical and the case $u=K$ is trivial.
\end{proof}
A coinvariant $\CB'_{\al}:\CP_{\al}\ot\CP_{\al}$ is then obtained by
combining $\CB_{\al}$ with the intertwining operator $\CI_{\al}$:
\begin{equation}
\CB_{\al}'\;\equiv\; \CB_{\al}\circ(\CI_{\al}\ot\id).
\end{equation}

In order to construct coinvariants $\CB^{(n)}$ for $n>2$ one may use 
intertwining maps
\[
\CC\in{\rm Hom}_{\USL}(\CP_{\al_{n-1}}\ot\ldots\ot\CP_{\al_1}, \CP_{\al_n}).
\]
Such maps can be constructed by iterating Clebsch-Gordan maps, as
has been discussed explicitly in the case $n=4$ at the beginning of the
present Section. One may associate a coinvariant $\CB_{\CC}$
to any $\CC\in{\rm Hom}_{\USL}(\CP_{\al_{n-1}}\ot\ldots\ot\CP_{\al_1}, 
\CP_{\al_n})$  via
\begin{equation}
\CB_{\CC}^{}\;\equiv\; \CB\circ(\id\ot\CC).
\end{equation}

The maps $\CC$ can be represented 
explicitly with the help of meromorphic integral kernels 
$\Phi_{\CC}(x_{n};\fx)$, $\fx\equiv(x_{n-1},\ldots,x_1)$ that generalize
$\Phi^{\flat}_{\al_{\flat}}$ and the Clebsch-Gordan coefficients. 
It follows that the corresponding coinvariant $\CB_{\CC}$ can 
be represented as 
\begin{equation}\label{coinvint}
\CB_{\CC}(f_n\ot\ldots\ot f_1)\;=\; 
\int\limits_{\BR}dx_n \;T_{x_n}^{i\frac{Q}{2}}f_n(x_n) \; 
\int\limits_{\BR^{n-1}}d\fx \;\Phi_{\CC}(x_{n};\fx) \;f_{n-1}(x_{n-1})\ldots
f_1(x_1).
\end{equation}
It is possible to rewrite \rf{coinvint} as a convolution of
$f_{n}(x_{n})\ldots
f_1(x_1)$ against a kernel $\Psi_{\CC}(\ux)$, $\ux\equiv(x_n,\ldots,x_1)$:
To this aim it is necessary to ``partially'' integrate the finite 
difference operator in \rf{coinvint} to let it act on $\Phi_{\CC}$.
One should note that the analytic continuation 
of the integral over $\fx$ to complex values of $x_n$ 
may in general be represented by integrating the variable $\fx$
over deformed contours, cf. e.g. the proof of Proposition \ref{intertwprop}.
One arrives at a representation of the form
\begin{equation}
\CB_{\CC}(f_n\ot\ldots\ot f_1)\;=\; 
\int\limits_{C^{n}}dx_n\ldots dx_1 
\;\Psi^{}_{\CC}(x_n,\ldots,x_1)f_{n}(x_{n})\ldots
f_1(x_1),
\end{equation}
where 
\begin{equation}\Psi^{}_{\CC}(x_n,\ldots,x_1)\;=\;
T_{x_n}^{-i\frac{Q}{2}}\;\Phi_{\CC}(x_{n};x_{n-1},\ldots,x_1).
\end{equation}
\begin{rem}
The kernels that represent the coinvariants are in some 
respects analogous to functional realizations of the conformal blocks 
in conformal field theory. We strongly suspect that we are 
touching upon the tip of an iceberg at this point: Quantization
of Teichm\"uller space, as developed in \cite{Fo}\cite{Ka}
conjecturally leads to a construction of spaces of conformal blocks
in Liouville theory. One may expect this to be equivalent to 
a quantization of certain moduli spaces of flat $SL(2,\BR)$
connections on Riemann surfaces with marked points. In analogy to
results of \cite{AS} one would expect spaces of conformal blocks 
in the case of the punctured Riemann sphere to be represented
by spaces of coinvariants in tensor products of $\USL$ representations.
A class of these has been constructed in the present subsection.
It would certainly be rather interesting and far-reaching if one
could establish a direct relation between these spaces and the 
Hilbert spaces constructed via quantization
of Teichm\"uller space. 

In this regard we find 
the following observation quite intriguing: 
Consider the case of $n=4$. There is a canonical way to define a
Hilbert space $\CH^{(0,4)}$ of coinvariants by taking 
the sets $\{\Phi_{\al}^{\flat};\al\in\BS\}$ for either
$\flat=s$ or $\flat=t$ as basis in the
sense of generalized functions with the normalization given by 
\begin{equation}
(\,\Phi_{\al}^{\flat}\,,\,\Phi_{\al'}^{\flat}\,)\;=\,
|S_b(2\al)|^{-2}\de(\al-\al').
\end{equation}  
The observation made in subsection 5.6. now implies
that $\CH^{(0,4)}$ is in a canonical way isomorphic
to $L^2(\BR)$ such that multiplication with
$[\al_{s}-\frac{Q}{2}]_b^2$ (resp. $[\al_{s}-\frac{Q}{2}]_b^2$)
gets mapped into the self-adjoint
finite difference operator $\CQ_{s}$ (resp. $\CQ_{t}$). Maybe 
there is a rather direct connection of these operators
to the geodesic length operators appearing in the quantization 
of Teichm\"uller space. This would establish a direct relation between
the latter and our quantum group results.
\end{rem}

\newpage

\section{Appendix A: Spectral analysis of $C_{21}(\ka_3)$}

This appendix is devoted to the proof of Theorem \ref{casspec}.

\subsection{Preliminaries}

The difference operator to be considered is of the form
\begin{equation}\label{Cas}
C_{21}(\ka_3)-[\al_3-\fr{Q}{2}]^2_b\;=\;
\de_{+}e^{\pi i bQ}e^{2\pi b x}-\de_0+\de_{-}e^{-\pi i bQ}
e^{-2\pi b x},
\end{equation}
where $\de_s$, $s=-,0,+$ are x-independent finite difference operators
given by
\begin{equation}\begin{aligned}
\de_+= & \;T_x^{-ib}{[}d_x-\al_2-ik_3{]}_b^{}
         {[}d_x-\al_1+ik_3{]}_b^{}\\
2\de_0\;= & \;\{0\}_b^{}\Bigl( \{Q\}_b^{}T_x^{-2ib}-
\bigl(e^{-2\pi bk_3}\{2\al_2-Q\}_b^{}+e^{2\pi bk_3}\{2\al_1-Q\}_b^{}\bigr)
T_x^{-ib}+\{2\al_3-Q\}_b^{}\Bigr)\\
\de_-= & \;T_x^{-ib}{[}d_x+\al_2-ik_3{]}_b^{}
{[}d_x+\al_1+ik_3{]}_b^{},
\end{aligned}\end{equation}
and $\ka_3=-2k_3$.
It will initially be defined on the domain $\FD\subset L^2(\BR)$ consisting
of functions with the following property: There exists a function $F(z)$ that 
is \begin{enumerate}
\item  holomorphic 
in the strip $\{ z\in \BC|\Im(z)\in[-2b,0]\}$ and
\item the functions $F_y(x)\equiv F(x+iy)$ are in $L^2(\BR,dx\cosh(2\pi bx))$
for any $y\in[-2b,0]$.
\end{enumerate}
\begin{propn}\label{symm}
The operator $(C_{21}(k_3),\FD)$ is a symmetric, densely defined operator in 
$L^2(\BR)$. The domain $\FD^{\dagger}$ of its adjoint is dense as well.
\end{propn}
\begin{proof} 
First of all note that one has 
\begin{equation}
(f,T_x^{-ib} g) = ( T_x^{-ib}f,g)
\end{equation}
for any $f,g\in\FD$. This follows by shifting  
the contour of the integration that represents $ (f,T_-g)$ to the 
line $\BR+ib$. The fact that $C_{21}(\ka_3)$ is symmetric is then seen
by a simple calculation remembering that $\al_i^*=Q-\al_i$, $i=1,2$.

The fact that $\FD$ and $\FD^{\dagger}$ are dense in $L^2(\BR)$
is easily seen by
noting that any Hermite-function is contained in these
sets.
\end{proof}

The Paley-Wiener theorem provides a characterization of the Fourier-transform
$\tilde{\FD}$ of the domain $\FD$ of $C_{21}(\ka_3)$. The action 
of $C_{21}(\ka_3)$ on functions in $\FD$ then corresponds to acting on 
$\tilde{\FD}$ with the following operator:
\begin{equation}\label{CasFT}\begin{aligned}
C_{21}(\ka_3) -\bigl[ &  \al_3  -\fr{Q}{2}  \bigr]_b^2\;\equiv \;
\De_0-e^{2\pi b \om}\De_1+
e^{4\pi b \om} \De_2\\[1ex]
\De_0 =
& \;{[}d_{\om}+\al_3-Q-\fr{1}{2}(\al_1+\al_2){]}_b^{}
                 {[}d_{\om}-\al_3-\fr{1}{2}(\al_1+\al_2){]}_b^{}\\
\De_1 = & \;{[}d_{\om}+\fr{1}{2}(\al_1+\al_2){]}_b^{}\Bigl(
e^{i\pi b(d_{\om}-\frac{1}{2}(\al_1+\al_2)+Q)}\{\al_1-\al_2-2ik\}_b^{}\\
& \qquad\qquad\qquad\qquad\qquad\quad
-e^{-i\pi b(d_{\om}-\frac{1}{2}(\al_1+\al_2)+Q)}\{\al_1-\al_2+2ik\}_b^{}
\Bigl)\\
\De_2 = & \;{[}d_{\om}+\fr{1}{2}(\al_1+\al_2){]}_b^{}
{[}d_{\om}+\fr{1}{2}(\al_1+\al_2)+Q{]}_b^{}.
\end{aligned}
\end{equation}

\subsection{Strategy}

The key to the proof of Theorem \ref{casspec} 
is the following result characterizing 
regularity and asymptotic properties of distributional solutions
to the eigenvalue equation of the operator $C_{21}(\ka_3)$:
\begin{thm}\label{regthm}
Let $\Phi\in\CS'(\BR)$ be a distributional solution of 
$(C_{21}(\ka_3)-[\al_3-\frac{Q}{2}]^2)^t\,\Phi=0$. 
\begin{enumerate}\item
$\tilde{\Phi}$ is represented by a function $\tilde{\Phi}(\om)$ that
can be continued to a meromorphic function on $\BC$, 
with simple
poles within $\FS_{Q/2}$ only at 
\[
\begin{aligned}
\om=&-k_3+i(\al_1+nb+mb^{-1}),\\
\om=&+k_3+i(\al_2+nb+mb^{-1}),
\end{aligned}
\quad
\begin{aligned}
\om=&-k_3-i(\al_1+nb+mb^{-1}),\\
\om=&+k_3-i(\al_2+nb+mb^{-1}),
\end{aligned}
\quad
n,m\in\BZ^{\geq 0}.
\]
\item $\Phi$ can be represented 
as $\Phi=\lim_{\ep\ra 0}\Phi_{\ep}$ where 
$\Phi_{\ep}$ is for 
$\ep>0$ represented as the restriction to $\BR$ of a function 
$\Phi_{\ep}(x)$ that is meromorphic on $\BC$
with poles only at 
\[ \begin{aligned}
x= &  
+\fr{i}{2}\bigl(\al_1+\al_2-Q\bigr)\pm i\bigl(\al_3-\fr{Q}{2}\bigr)
-i(\ep+nb+mb^{-1}),\\
x=&-\fr{i}{2}\bigl(\al_1+\al_2-Q\bigr)+i\bigl(\fr{Q}{2}+nb+mb^{-1}\bigr),
\end{aligned}\qquad
n,m\in\BZ^{\geq 0}.
\]
\end{enumerate}
\end{thm}

In fact, given these properties it is not very difficult to show that
for any given eigenvalue $[\al_3-\frac{Q}{2}]^2$ 
there is at most one tempered distributional solution 
to the eigenvalue equation (Proposition \ref{uniprop}). Moreover, no such 
solution exists for $\Re(2\al_3-Q)\neq 0$. It follows \cite{AG} that the
deficiency indices vanish and 
$C_{21}(\ka_3)$ has a unique self-adjoint extension. The spectral decomposition 
can be written as expansion into generalized eigenfunctions \cite{GV}.
It can be shown on rather general grounds that only tempered distributions
can appear in the spectral decomposition, as nicely discussed in \cite{Be}.
The combination of Theorem \ref{regthm} and Proposition \ref{uniprop}
therefore also yields a characterization of the support of the 
Plancherel measure.

These remarks reduce the proof of Theorem \ref{casspec} to that of Theorem \ref{regthm}
and Proposition \ref{uniprop}.

\subsection{Preparations}

In view of the explicit expressions for $C_{21}(\ka_3)$ (cf. \rf{Cas}) resp.
its Fourier-transform \rf{CasFT} one may anticipate
that the analysis of the asymptotic behavior of $\Phi$ and 
$\TP$ will require some information about properties of the
operators $\de_+$, $\de_-$ resp. $\De_0$, $\De_2$. 
The information that will be needed is contained in the 
following Lemmas:

\begin{lem}\label{deplinv} $\de_{\pm}$ is invertible on 
$\CC_c^{\infty}(\BR)$. The image $f(x)$ of a 
function $g\in\CC_c^{\infty}(\BR)$ under $\de_{\pm}^{-1}$ 
has the following properties:
\begin{enumerate}
\item $f(x)$ is analytic in the strip $\{x\in\BC;\Im(x)\in(-2b,0)\}$ and
$f(x)\in\CC^{\infty}(\BR)$,  $f(x-2ib)\in\CC^{\infty}(\BR)$.
\item $\tilde{f}(\om)$ is meromorphic in $\BC$ with simple poles at
\[
\om=-k_3+i(\mp\al_1+nb^{-1})\qquad
\om=+k_3+i(\mp\al_2+nb^{-1})\qquad n\in\BZ.
\]
\end{enumerate}\end{lem}
\begin{proof}
The action of $\de_{\pm}^{-1}$ is represented on the Fourier transform
$\tilde{f}$ as multiplication with 
\[
(\tilde{\de}_{\pm})^{-1}(\om)\equiv e^{-2\pi b\om}
{[}i\om\mp\al_2-ik_3{]}_b^{-1}
{[}i\om\mp\al_1+ik_3{]}_b^{-1} . 
\]
The statement on the 
analyticity properties of $\tilde{f}$ is then clear  
after recalling that
the function 
$\tilde{g}(\om)$ is entire analytic and of rapid decay being the
Fourier transform of a $\CC_c^{\infty}$ function \cite[Theorem IX.11]{RS2}. 

The statement that $(\de_+^{-1}g)(x)$ is analytic in the strip 
$\{x\in\BC;\Im(x)\in(-2b,0)\}$ follows from the asymptotic decay properties
of $(\tilde{\de}_{\pm}^{-1})(\om)$  by means of the Paley-Wiener Theorem. 
In fact, the rapid decay of $\tilde{g}(\om)$ ensures convergence
of the inverse Fourier transformation for any
x-derivative of $(\de_+^{-1}g)(x)$ even in the extremal cases
$\Im(x)=0$ and $\Im(x)=-2b$.
\end{proof}

We will furthermore need similar statements about the inverses of
$\De_0$ and $\De_2$. 
\begin{lem} \label{Deinv} $\De_2$ is invertible on 
$\CC_c^{\infty}(\BR)$. The image $f(\om)$ of a 
function $g\in\CC_c^{\infty}(\BR)$ under $\De_{2}^{-1}$ 
has the following properties:
\begin{enumerate}
\item $\tilde{f}(x)$ is meromorphic in $\BC$ with simple poles at
\[ x=-\fr{i}{2}(\al_1+\al_2)-i(Q+nb^{-1})
\qquad 
x=-\fr{i}{2}(\al_1+\al_2)+inb^{-1}\qquad n\in\BZ.
\]
\item $f(\om)$ is 
analytic in the strip $\{\om\in\BC;\Im(x)\in(-b,b)\}$ and
$f(\om\pm ib)\in\CC^{\infty}(\BR)$.
\end{enumerate}
\end{lem}
\begin{lem} $\De_0$ is invertible on the space of functions
\[ 
\cD(\De_0)\equiv \bigl(d_{\om}+\al_3-Q-\fr{1}{2}(\al_1+\al_2)\bigr)
\bigl(d_{\om}-\al_3-\fr{1}{2}(\al_1+\al_2\bigl)h, \quad 
h\in\CC_c^{\infty}(\BR). 
\]
The image $f(\om)$ of a 
function $g\in\cD(\De_0)$ under $\De_{0}^{-1}$ 
has the following properties:
\begin{enumerate}
\item $\tilde{f}(x)$ is meromorphic in $\BC$ with simple poles at
\[ x=+\fr{i}{2}(\al_1+\al_2-Q)\pm i(\al_3-\fr{Q}{2})-inb^{-1}
\qquad n\in\BZ\setminus {\{} 0 {\}}.
\]
\item $f(\om)$ is 
analytic in the strip $\{\om\in\BC;\Im(x)\in(-b,b)\}$ and
$f(\om\pm ib)\in\CC^{\infty}(\BR)$.
\end{enumerate}
\end{lem}

\subsection{Asymptotic estimates}

We now want to show that the Fourier-transform $\TP$ of $\Phi$ may 
actually be represented by integration against a function $\TP(\om)$.
For technical reasons it will be necessary to start by considering the 
distribution $\Phi_R\in\CS'(\BR)$ defined by 
\[
\TP_R^{} \equiv \tilde{\de}_{\rm tr,R}(\om)\TP
\equiv \prod_{\substack{\om'\in\CI_+\cup\CI_-\\
|\Im(\om')|<R}}(\om-\om')\;\TP,
\]
where $\CI_+$ (resp. $\CI_-$) are the sets of values for $\om$ where either
$\tilde{\de}_+(\om)$ 
or $\tilde{\de}_-(\om)$ have a pole in the upper (resp. lower) half plane. 
The following result characterizes the asymptotic behavior of $\Phi_R$.
\begin{propn}\label{asest1} Let $\tau_n\in\CC_c^{\infty}(\BR)$ 
have support only in $[n-1,n+1]$. For sufficiently large value of $R$ 
there exists some $N>0$ such that 
\begin{equation}
\cosh(2\pi b n)\bra \Phi_R,\tau_n\ket < N \quad \text{for all $n\in\BZ$}.
\end{equation}
\end{propn}
\begin{proof}
We will rewrite $\bra \Phi_R,\tau_n\ket$ in a form that allows us to estimate
its asymptotics for large $n$. One may write
\begin{equation}\label{as1}\begin{aligned}
\bra\Phi_R,\tau_n\ket=& \bra\Phi,\de_{\rm tr,R}\tau_n\ket,\\
=& \bra\Phi,\de_+e^{2\pi bx} \si_{n,R}\ket, \qquad\quad \text{where} 
\quad \si_{n,R}
\equiv e^{-2\pi bx}(\de_+)^{-1}\de_{\rm tr,R}\tau_n;\\
=& \bra\Phi,\de_+^c\si_{n,R}\ket,\qquad \qquad \quad \text{where} \quad 
\de_+^c\equiv (\de_0-\de_-e^{-2\pi b x}).
\end{aligned}
\end{equation}
In the last step we have used that $\Phi$ weakly solves
the eigenvalue equation, for which one needs to check
that $\si_{n,R}\in\FD$: One point of having introduced $\de_{\rm tr,R}$ is that
it improves the asymptotic behavior of $(\de_+)^{-1}\de_{\rm tr,R}\tau_n$
for $x\ra -\infty$ by cancelling the poles of its Fourier transform
in $\{ \om\in\BC;\Im(\om)<R\}$.
 
The regularity theorem for tempered distributions 
\cite[Theorem V.10]{RS} allows us to 
furthermore write 
\begin{equation}\label{est2}
\bra\Phi_R,\tau_n\ket=\int\limits_{-\infty}^{\infty}dx \;\Theta(x)\;\rho_{n,R}(x)
\qquad \text{where} \quad \rho_{n,R}\equiv 
\pa_x^k \;\de_+^c e^{-2\pi bx}(\de_+^{})^{-1}\de_{\rm tr,R}\tau_n.
\end{equation}
for some positive integer $k$ and a polynomially bounded continuous 
function $\Theta(x)$. The functions $\rho_{n,R}(x)$ may be represented by  
expressions of the form
\begin{equation}
\rho_{n,R}(x)\;= \;  \sum_{k=1,2} \;\;
C_k e^{-2\pi bx} \int\limits_{-\infty}^{\infty}
d\om\;e^{2\pi i \om x}
\frac{P_{k,R}(\om)\tilde{\tau}_n(\om)}{(1-e^{2\pi b (\om-k+i\al_1)})
(1-e^{2\pi b (\om+i\al_2)})},
\end{equation}
where $P_{k,R}(\om)$ $k=1,2$ are some polynomials in $\om$. The 
functions $\rho_{n,R}(x)$ have main support around $x=n$, and by choosing
$R$ large enough one can achieve
decay stronger than $e^{-2\pi \la |x-n|}$ for any $\la>0$.  
It is then convenient to split the integral in \rf{est2} into an
integral $J_n$ obtained by integrating over $[\frac{n}{2},\frac{3n}{2}]$ 
and the remainder $J_n^c$. 

In order to estimate $J_n^c$ one may
use the polynomial boundedness of $\Theta(x)$ to estimate its 
absolute value by some constant times $\cosh(\ep x)$, where $\ep$ can 
be as small as one likes. The absolute value of $\rho_{n,R}(x)$ can in 
$\BR\setminus[\frac{n}{2},\frac{3n}{2}]$ 
be estimated by some inverse power of 
$\cosh(x)$, which is bounded by the chosen value of $R$. It follows that
the exist $D_1$, $N_1$ such that 
\begin{equation}
 |J_n^c|\;\leq\; D_1 e^{-2\pi \mu n} \quad \text{for any $n>N_1$}, 
\end{equation}
where $\mu$ can be made arbitrarily large by choosing $R$ large enough.

In the case of $J_n$ one may estimate $|\rho_{n,R}(x)|$ by
some constant times $e^{-2\pi b n}e^{-2\pi b|x-n|}$ and $\Theta(x)$
simply by a constant, which easily gives existence of $D_2$, $N_2$
such that  
\begin{equation}
\label{dest1} 
|J_n|\;\leq\; D_2 e^{-2\pi b n} \quad \text{for any $n>N_1$}. \end{equation}
This proves the claim about the asymptotics for $n\ra\infty$.
In the case of $n\ra-\infty$ one uses the operator
$\de_-$ in a completely analogous fashion
\end{proof}

\subsection{Representation of $\TP$} \label{TPrepr}

Assume that the set $\{\tau_n;n\in\BZ\}$ represents a 
$\CC_c^{\infty}(\BR)$-partition of 
unity. It will be convenient to choose the $\tau_n$ as translates of
$\tau_0$: $\tau_n(x)=\tau_0(x-n)$. This can always be achieved:
Let 
\begin{equation}\begin{aligned}
{} & \qquad\qquad\qquad\quad
\tau_0(x)=\left\{ \begin{aligned} {}& 0 & \quad\text{if}& \quad|x|>\fr{3}{4}\\
 & 1 & \text{if}& \quad |x|<\fr{1}{4}\\
 & \chi(x+\fr{1}{2}) & \text{if}& \quad x\in {[}-\fr{3}{4}, -\fr{1}{4}{]}\\
1- & \chi(x-\fr{1}{2})\; & \text{if}& \quad x\in {[}+\fr{1}{4}, +\fr{3}{4}{]},
\end{aligned}\right. \\
& 
\quad \chi(x)\;=\;N^{-1}\int\limits_{-\frac{1}{4}}^{x}dt \;\exp\biggl(\frac{1}{(x-\frac{1}{4})(x+\frac{1}{4})}\biggr)
\qquad N= \int\limits_{-\frac{1}{4}}^{\frac{1}{4}}dt \;\exp\biggl(\frac{1}{(x-\frac{1}{4})(x+\frac{1}{4})}\biggr)\end{aligned}
\end{equation}
The result of Proposition \ref{asest1} implies convergence of the following
sum
\begin{equation}\label{THdef}
\TP_R(\om)\;\equiv \;\sum_{n\in\BZ}\;\bra \Phi_R,\tau_n e^{-2\pi i\om x}\ket
\end{equation}
which defines $\TP_R(\om)$ as a function that is analytic in the strip
$\{\om\in\BC;\Im(\om)\in(-b,b)\}$.

\begin{propn} \label{TPprop}
The function $\TP_R(\om)$ represents the distribution $\Phi_R$
in the sense that 
\begin{equation} \label{Phirepr}
\bra\Phi_R,f\ket\;=\;\int\limits_{-\infty}^{\infty}
d\om\;\TP_R(\om)\tilde{f}(\om).
\end{equation}
\end{propn}
\begin{proof} 
To begin with, note that $\Phi_{R,n}(\om)\equiv 
\bra \Phi_R,\tau_n e^{-2\pi i\om x}\ket$ represents the Fourier-transform
of the distribution $\tau_n \Phi_R\in\CS'(\BR)$ of compact support 
\cite[Theorem IX.12]{RS2}.  It follows that 
$\bra\Phi_R,\tau_n e^{-2\pi i\om x}\ket$
is polynomially bounded. Since the convergence in \rf{THdef} is absolute,
one concludes that $\TP_R(\om)$ is polynomially bounded as well. In the
evaluation of $\TP_R(\om)$ against a test-function $f\in\CS(\BR)$
one may therefore insert definition \rf{Phirepr} and exchange the orders
of integration and summation to get
\begin{equation}\begin{aligned}
\int\limits_{-\infty}^{\infty}d\om\;\TP_R(\om)\tilde{f}(\om)  
\;=& \; \sum_{n\in\BZ} \;\,\int\limits_{-\infty}^{\infty} 
d\om\; \TP_{R,n}(\om)\tilde{f}(\om) \\
\;=& \; 
\sum_{n\in\BZ}\;\bra \Phi_R,\tau_n f\ket
\;=\; \bra\Phi_R,f\ket,
\end{aligned}\end{equation}
where we used that fact that the set $\{\tau_n;n\in\BZ\}$ represents a 
partition of unity in the last step.
\end{proof}

In order to recover the sought-for distribution $\Phi$ from $\Phi_R$
one only has to divide $\TP_R(\om)$ by $\tilde{\de}_{\rm tr,R}(\om)$.
The resulting function is {\it meromorphic} in the strip 
$\{\om\in\BC;\Im(\om)\in(-b,b)\}$, with poles at distance
$\frac{1}{2}(b^{-1}-b)$ from the real axis. 

\subsection{Representation of $\Phi$}

In order to get a similar result on the representation of $\Phi$ in 
x-space we will analogously consider the asymptotics of $\TP$ in 
$\om$-space. Here it will be convenient to start by considering
\[
\Phi'_R \equiv \tilde{\de}'_{\rm tr,R}(x)\Phi
\equiv \prod_{s\in\{+,-\}}(x-x_{s})\prod_{\substack{y\in\CI_+\cup\CI_-\\
|\Im(z)|<R}}(x-y)\;\Phi,
\]
where $\CI_+$ (resp. $\CI_-$) 
denotes the union of the sets of zeros of $\tilde{\De}_2(z)$ and
$\tilde{\De}_0(z)$ which lie in the 
upper (resp. lower) half plane,
and $x_{\pm}$ are the zeros of $\tilde{\De}_0(z)$ that lie {\it on} 
the real axis, given by
\[ x_{\pm}\equiv
+\fr{i}{2}\bigl(\al_1+\al_2-Q\bigr)\pm i\bigl(\al_3-\fr{Q}{2}\bigr).
\]

For the asymptotics of $\TP_R'$ one has a result completely analogous to 
Proposition \ref{asest1}:
\begin{propn}\label{asest2}
Let $\{\tau_n;n\in\BZ\} $ be a sequence of functions in
$\CC_c^{\infty}(\BR)$ that have support only in $[n-1,n+1]$.
For sufficiently large $R$ there exists some $N>0$ such that 
\begin{equation}
\cosh(2\pi b n)\bra \TP'_R,\tau_n\ket < N \quad \text{for all $n\in\BZ$}.
\end{equation}
\end{propn}
\begin{proof} The proof is to a large extend analgous to that of
Proposition \ref{asest1}, so we will only sketch some necessary 
modifications. 

In order to get an estimate of $\bra \TP'_R,\tau_n\ket$ for $n\ra -\infty$
one may use the eigenvalue equation to rewrite it as
\begin{equation}\begin{aligned}
\bra\TP'_R,\tau_n\ket= &\bra \TP,
\De_0^{}\De_0^{-1}\de'_{\rm tr,R}\tau_n\ket\\
=&\bra \TP,\De_0^{c} \De_0^{-1}\de'_{\rm tr,R}\tau_n\ket \qquad
\text{where $\De_0^{c}=e^{2\pi b\om}\De_1-e^{4\pi b \om}\De_2$}.
\end{aligned}\end{equation}
It follows as in the proof of Proposition \ref{asest1} that
$\bra\TP'_R,\tau_n\ket\sim e^{+2\pi bn}$ for $n\ra -\infty$. 

In the case of $n\ra\infty$ one may use instead
\begin{equation}\begin{aligned}
\bra\TP'_R,\tau_n\ket= 
&\bra \TP,e^{4\pi b\om}\De_2^{}\De_2^{-1}e^{-4\pi b\om}\de'_{\rm tr,R}
\tau_n\ket\\
= & \bra \TP,\De_2^{c} 
\De_2^{-1}e^{-4\pi b\om}\de'_{\rm tr,R}\tau_n\ket \qquad
\text{where $\De_0^{c}=e^{2\pi b\om}\De_1-\De_0$},
\end{aligned}\end{equation}
which gives $\bra\TP'_R,\tau_n\ket
\sim e^{-2\pi bn}$ for $n\ra \infty$. 
\end{proof}

It follows as in the previous
section that $\Phi_R'$ is represented by convolution against
a function $\Phi_R'(x)$ which is holomorphic in 
$\{x\in\BC;\Im(x)\in(-b,b)\}$. In this case, however, recovering
$\Phi$ from $\Phi_R'$ is more subtle since $\tilde{\de}'_{\rm tr,R}(x)$ has
two simple zeros on the real axis. The resulting ambiguity in the 
definition of $\Phi$ in terms of $\Phi_R'(x)$ is well-known
(cf. e.g. 
\cite[Chapter V, Example 9]{RS}) and may be parametrized as follows:
\begin{equation}\label{gendiv}
\Phi=\prod_{s\in\{+,-\}}\biggl(\frac{C_{s}}{x-x_{s}+i0}+
\frac{1-C_{s}}{x-x_{s}-i0}\biggr)
\prod_{\substack{y\in\CI_+\cup\CI_-\\
|\Im(z)|<R}}\frac{1}{x-y}\;\Phi_R'(x).
\end{equation}

Lemma \ref{distPW} then describes the corresponding asymptotic
behavior of $\TP(\om)$. 
In general one would find terms 
with exponential decay weaker than $e^{-2\pi b |\om|}$ for
$\om\ra\infty$ that come either
from zeros of $\tilde{\de}'_{\rm tr,R}(x)$ strictly
above the real axis, or from $x_{\pm}$ in the case of
$C_s\neq 0$. The occurrence of such terms
can be excluded by means of the following argument: 
\begin{lem}\label{vanres}
Let $\Phi\in\CS'(\BR)$ be a distributional solution of 
$(C_{21}(\ka_3)-[\al_3-\frac{Q}{2}]^2)^t\,\Phi=0$ that 
is represented by a function $\TP(\om)$ which has asymptotic behavior
for $\om\ra\infty$ of the form 
\[
\TP(\om)\;=\;  +2\pi i \sum_{j\in\CI_{-}}e^{-2\pi iz_j\om}R_j+\TP_{a_-}(\om),
\]
where $\TP_{b}(\om)$ decays at least as fast as  $e^{-2\pi b\om}$ for
$\om\ra \infty$. Then $R_j=0$ if $\Im(z_j)<b$.
\end{lem}
\begin{proof} 
Consider $\bra \TP,\tau_n\ket$, where now $\tau_n$ is chosen proportional
to $e^{-\ka(x-n)^2}$. One has
\begin{equation}\label{weakev}
\bigl[   \al_3  -\fr{Q}{2}  \bigr]_b^2\;\bigl\bra \TP,\tau_n\bigr\ket
\;=\;\Bigl\bra \TP, \Bigl(\De_0-e^{2\pi b \om}\De_1+
e^{4\pi b \om}\De_2+\bigl[   \al_3  -\fr{Q}{2}  \bigr]_b^2\Bigr)\tau_n
\Bigr\ket.
\end{equation}
Now if there were terms
with exponential decay weaker than $e^{-2\pi b \om}$ in the asymptotic
expansion of $\TP(\om)$ for $\om\ra \infty$ one would find
terms terms that grow exponentially with $n\ra \infty$ on the 
right hand side of \rf{weakev}. But polynomial boundedness of 
$\TP$ excludes the occurrence of such terms on the left hand side
of \rf{weakev}. 
\end{proof} 

\subsection{Completing the proof of Theorem \ref{regthm}}

Concerning the distribution $\Phi$, we previously found that
away from its singular support at $x=x_{\pm}$ it is 
represented by a function $\Phi(x)$. The asymptotic behavior 
of $\Phi(x)$ is via Lemma \ref{distPW} given by the analytic 
properties of $\TP$ that were stated after the proof of Proposition
\ref{TPprop}. 
The possible
poles of $\TP$ at distance $\frac{1}{2}(b^{-1}-b)$ from the real axis
would lead to terms which decay more slowly as $e^{-2\pi b|x|}$ for
$|x|\ra\infty$.  
The appearance of such terms
can now easily be excluded by an argument analogous to
the proof of Lemma \ref{vanres} in the x-representation.

Furthermore, knowing that the function $\Phi(x)$ that represents
$\Phi$ away from its singular support
decays exponentially for
$|x|\ra\infty$ allows us to use an argument very similar to the
proof of Proposition \ref{asest1} to further improve upon the
estimate of the rate of decay as given in Proposition \ref{asest1}:
In estimating $J_n$ one may for large enough $n$ replace $\Theta(x)$
by $\Phi(x)$. The exponential decay of the latter may then be used to 
improve \rf{dest1} to 
\begin{equation}
\label{dest2} 
|J_n|\;\leq\; D_2 e^{-2\pi \nu n} \quad \text{for any $n>N_1$}. \end{equation}
for some $\nu > b$, implying that
$\Phi(x)$ decays faster than $e^{-2\pi b|x|}$ for
$|x|\ra\infty$.

But this means via Lemma \ref{distPW} that the Fourier-transformation
$\TP(\om)$ is analytic in an open strip containing 
$\{\om\in\BC;|\Im(\om)|<b\}$, and that $\TP(\om)$ solves 
$(\tilde{C}_{21}(k_3) -[  \al_3  -\fr{Q}{2}  ]_b^2)^t
\TP(\om)=0$ in the ordinary sense. The meromorphic extension to all of 
$\BC$ is then easily obtained by using the eigenvalue equation to define
the values of $\TP(\om)$ outside $\{\om\in\BC;|\Im(\om)|<b\}$ in terms
of those inside. This finishes the proof of the first half of 
Theorem \ref{regthm}. The completion of the proof of the second 
half proceeds along very similar lines.

\subsection{Uniqueness of generalized eigenfunctions}

Theorem 3 also implies that the meromorphic function 
$\Phi(x)$ that represents
the distribution $\Phi$ 
must solve the transpose of the eigenvalue equation
in the usual sense.

\begin{propn}\label{uniprop} 
There is at most one solution to 
$(C_{21}(\ka_3)-[\al_3-\frac{Q}{2}]^2)^t\,\Phi(x)=0$
that has the analytic and asymptotic properties that follow from 
Theorem \ref{regthm}.
\end{propn}
\begin{proof}
If one introduces $\Xi(x)$ via 
(recall $\ka_3=-2k_3$)
\begin{equation}\begin{aligned}
\Phi(x)=e^{\pi x (\al_3+\al_1-\al_2-i\ka_3)}
& \frac{S_b(-ix-\frac{1}{2}(\al_1+\al_2)+\al_3)}{S_b(-ix+\frac{1}{2}
(\al_1+\al_2))}\\
&\quad \qquad\quad \qquad
\ti \Xi\bigl(x-\fr{i}{2}(\al_1+\al_2-2(Q-\al_3))\bigr),
\end{aligned}\end{equation}
one may verify by direct calculation using the functional equation
of the function $S_b(x)$ that the equation 
$(C_{21}(\ka_3)-[\al_3-\frac{Q}{2}]^2)^t\,\Phi(x)=0$ is equivalent to the 
following equation for $\Xi(x)$:
\begin{equation}\label{redeq}\begin{aligned}
\Bigl( (1-e^{2\pi i b(\al_3+\al_1-\al_2)}T_x^{ib}) & (1-e^{2\pi i b
(\al_3-i\ka_3)}T_x^{ib})\\
 & -e^{-2\pi bx}(1-T_x^{ib})(1-e^{2\pi i b(\al_1-\al_2- i\ka_3)}
T_x^{ib})\Bigr)\Xi(x)=0.\end{aligned}
\end{equation}
By using Lemma \ref{distPW} and the 
properties of $S_b(x)$ that are summarized in 
Appendix B one may deduce the following properties of the 
Fourier transform $\TX(\om)$ 
of $\Xi(x)$
from Theorem \ref{regthm}:
\begin{enumerate}
\item $\Xi(x)$ has a Fourier transform $\TX(\om)$ that 
is analytic
in $\{ \om\in\BC;{\rm Im}(\om)\in(-Q/2,0) \}$, and 
\item $\TX(\om)$ has the following asymptotic behavior for 
$\om\ra\pm \infty$:
\[
\TX(\om)= R_+(\om), \qquad
\TX(\om)= K_- + R_-(\om),
\]
where $K_-$ is a constant, $R_-(\om)$ has exponential decay for 
$\om\ra -\infty$ and
$R_+(\om)$ has exponential decay stronger than $e^{-4\pi b\om}$ 
for $\om\ra \infty$.
\end{enumerate}
Equation \rf{redeq} is equivalent to the following {\it first order}
difference equation for $\TX(\om)$:
\begin{equation}\label{redeqFT}\begin{aligned}
\Bigl( (1- & e^{2\pi i b(\al_3+\al_1-\al_2-i\om)})  
(1-e^{2\pi i b(\al_3-i\ka_3-i\om)})\\
  & \qquad -(1-e^{2\pi ib(Q-i\om)})(1-e^{2\pi i b(Q+\al_1-\al_2-i\ka_3-i\om)}
)T_{\om}^{ib}\Bigr)\TX(\om)=0.\end{aligned}
\end{equation}
Now there exists a solution to \rf{redeqFT}, namely
\begin{equation}
\TX(\om)=\frac{G_b(\al_3+\al_1-\al_2-i\om)G_b(\al_3-i\ka_3-i\om)}{G_b(
Q-i\om)G_b(Q+\al_1-\al_2-i\ka_3-i\om)}, 
\end{equation}
that has all the required analytic and asymptotic properties.
If there was a second solution $\TX'(\om)$ of these conditions one could
consider the ratio $Q(\om)\equiv \TX'(\om)/\TX(\om)$. This ratio 
must be a solution to $(T_{\om}^{ib}-1)Q(\om)=0$. Since 
$\TX(\om)$ has no zeros in the open strip
$\{ \om \in\BC;{\rm Im}(\om)\in(-Q/2,0) \}$ one concludes that
$Q(\om)$ is {\it holomorphic} in any such strip. 
The function $Q(\om)$ must furthermore be asymptotic to the 
constant function for $\om\ra\pm\infty$. 
But this implies that $Q={\rm const.}$: The 
function $P(z)\equiv Q(\frac{b}{2\pi}\ln(z))$ is holomorphic and regular
on the whole Riemann sphere, therefore constant.
\end{proof}
\newpage
\section{Appendix B: Special functions}

The basic building block for the class of special functions to be considered
is the Double Gamma function introduced by Barnes \cite{Ba},
see also \cite{Sh}.
The Double Gamma function is defined as
\begin{equation}
\log\Ga_2(s|\om_1,\om_2)=  \Biggl(\frac{\pa}{\pa t}\sum_{n_1,n_2=0}^{\infty}
(s+n_1\om_1+n_2\om_2)^{-t}\Biggr)_{t=0}.
\end{equation}
Let $\Ga_b(x)=\Ga_2(x|b,b^{-1})$, and define 
the Double Sine function $S_b(x)$ and the Upsilon function $\up_b(x)$
respectively by
\begin{equation}
S_b(x)=\frac{\Ga_b(x)}{\Ga_b(Q-x)}\qquad \up_b(x)=
\frac{1}{\Ga_b(x)\Ga_b(Q-x)}.
\end{equation}
It will also be useful to introduce
\begin{equation}
G_b(x)=e^{\frac{\pi i}{2}x(x-Q)}S_b(x).
\end{equation}

\subsection{Useful properties of $S_b$, $G_b$}
\subsubsection{Self-duality}
\begin{equation}
S_b(x)=S_{b^{-1}}(x) \qquad G_b(x)=G_{b^{-1}}(x).
\end{equation}

\subsubsection{Functional equations}
\begin{equation}\label{S:funrel}
S_b(x+b)=2\sin(\pi b x)S_b(x) \qquad G_b(x+b)=(1-e^{2\pi i bx})G_b(x).
\end{equation}

\subsubsection{Reflection property}
\begin{equation}\label{reflprop}
S_b(x)S_b(Q-x)=1\qquad G_b(x)G_b(Q-x)=e^{\pi i (x^2-xQ)}.
\end{equation}

\subsubsection{Analyticity}
$S_b(x)$ and $G_b(x)$ are meromorphic functions with poles at
$x=-nb-mb^{-1}$ and zeros at $x=Q+nb+mb^{-1}$, $n,m\in\BZ^{\geq 0}$.

\subsubsection{Asymptotic behavior}
\begin{equation}\label{SGas}
S_b(x) \sim 
\left\{ 
\begin{aligned} 
{} & e^{-\frac{\pi i}{2}(x^2-xQ)} & 
\text{ for }\Im(x)\ra +\infty \\
   & e^{+\frac{\pi i}{2}(x^2-xQ)} & 
\text{ for }\Im(x)\ra -\infty
\end{aligned}
\right.\qquad
G_b(x) \sim 
\left\{ 
\begin{aligned} 
{} & 1 & 
\text{ for }\Im(x)\ra +\infty \\
   & e^{+\pi i(x^2-xQ)} & 
\text{ for }\Im(x)\ra -\infty
\end{aligned}
\right.
\end{equation}

\subsection{b-beta integral}

\begin{lem}
We have
\begin{equation}\label{bbeta}
B_b(\al,\be)\;\equiv\;\frac{1}{i}
\int\limits_{-i\infty}^{i\infty}d\tau \;e^{2\pi i \tau \be}
\frac{G_b(\tau+\al)}{G_b(\tau+Q)}\;=\;\frac{G_b(\al)G_b(\be)}{G_b(\al+\be)}
\end{equation}\end{lem}
\begin{proof}
From the relation (recall $T_{\tau}f(\tau)\equiv f(\tau+b)$)
\begin{equation}
0\;=\;\int\limits_{-i\infty}^{i\infty}d\tau (1-T_{\tau}^b)\;e^{2\pi i \tau \be}
\frac{G_b(\tau+\al)}{G_b(\tau+Q)},
\end{equation}
which easily follows from the analyticity and asymptotic properties 
of the $G_b$-function by means of Cauchy's theorem one 
finds the following functional equation for $ B_b(\al,\be)$:
\begin{equation}\label{funrelbb}
\frac{B_b(\al,\be+b)}{B_b(\al+b,\be)}=
\frac{1-e^{2\pi i b\be}}{1-e^{2\pi i b\be}}.
\end{equation}
By the $b\ra b^{-1}$ self-duality of $B_b$ one also has the same equation
with $b\ra b^{-1}$. For irrational values of $b$ it follows that 
\rf{funrelbb} and its $b\ra b^{-1}$ counterpart determine $B_b$ uniquely
up to a function of $\al+\be$. The expression  on the left hand side
of course satisfies \rf{funrelbb}. To fix the remaining ambiguity
one may note that the integral defining $B_b$ can be 
evaluated in the special case of $\al=b^{-1}$ by means of 
\cite[Chapt. 1.5., eqn. (28)]{E}:
\begin{equation}
B_b(b^{-1},\be)\;=\;\frac{b^{-1}}{1-e^{2\pi i b^{-1}\be}}.
\end{equation}
The equation \rf{bbeta} follows.
\end{proof}

Let us also introduce the combination 
\begin{equation}
\Theta_b(y;\al)\;\equiv\;\frac{G_b(y)}{G_b(y+\al)}.
\end{equation}
The b-beta-integral \rf{bbeta} can be read as a formula for the 
Fourier-transform of $\Theta_b(y;\al)$:
\begin{equation}
\Theta_b(y;\al)\;=\; \frac{1}{G_b(y)}\frac{1}{i}
\int\limits_{-i\infty}^{i\infty}d\tau\; e^{2\pi i \al\tau}\Theta_b(\tau+y;Q+y).
\end{equation}
An expansion describing the asymptotic behavior of $\Theta_b(y;\al)$
for $|\Im(y)|\ra\infty$ can therefore easily be obtained from Lemma 
\rf{distPW}: One finds
\begin{equation}\label{Thetaas}\begin{aligned}
\Theta_b(y;\al)\;\underset{\Im(y)\ra+\infty}{\simeq}\;&  
\sum_{n,m\geq 0}\Theta_{b,+}^{(n,m)}(\al)e^{2\pi i(nb+mb^{-1})y}\\
\Theta_b(y;\al)\;\underset{\Im(y)\ra-\infty}{\simeq}\;&  
\sum_{n,m\geq 0}\Theta_{b,-}^{(n,m)}(\al)e^{-2\pi i(\al+nb+mb^{-1})y},
\end{aligned}\end{equation}
where $\Theta_{b,+}^{(0,0)}(\al)=1$, $\Theta_{b,-}^{(0,0)}(\al)=
e^{-\pi i \al(\al-Q)}$. 

\subsection{b-hypergeometric function}

The b-hypergeometric function will be defined by an integral representation
that resembles the Barnes integral for the ordinary hypergeometric function:
\begin{equation}\label{Barnesint}
F_b(\al,\be;\ga;y)=\frac{1}{i}\frac{S_b(\ga)}{S_b(\al)S_b(\be)}
\int\limits_{-i\infty}^{i\infty}ds\;\, e^{2\pi i sy}
\frac{S_b(\al+s)S_b(\be+s)}{S_b(\ga+s)
S_b(Q+s) },
\end{equation}
where the contour is to the right of the poles at $s=-\al-nb-mb^{-1}$ and
$s=-\be-nb-mb^{-1}$ and to the left of the poles at $s=nb+mb^{-1}$
and $s=Q-\ga+nb+mb^{-1}$, $n,m=0,1,2,\ldots$.
The function $F_b(\al,\be;\ga;-ix)$ is a solution of the $q$-hypergeometric
difference equation
\begin{equation}
\bigl([\de_x+\al][\de_x+\be]-e^{-2\pi b x}[\de_x][\de_x+\ga-Q]\bigr)
F_b(\al,\be;\ga;-ix)=0, \qquad \de_x=\fr{1}{2\pi}\pa_x
\end{equation}
This definition of a q-hypergeometric function is closely related to
the one first given in \cite{NU}. 
\begin{lem}
Consider the case that $\Re(\al)=\Re(\be)=Q/2$, $\Re(\ga)=Q$.
$F_b(\al,\be;\ga;y)$ is analytic in $y$ in the strip 
$\{y\in\BC;\Re(y)\in(-Q/2,Q/2)\}$. 
The leading asymptotic behavior for $|\Im(y)|\ra\infty$ is given by
\begin{equation}
\begin{aligned}
F_b(\al,\be;\ga;y)=& 1+\CO(e^{2\pi iby}) +\\
 +& e^{2\pi i (Q-\ga) y}
\frac{S_b(\ga)}{S_b(2Q-\ga)}
\frac{S_b(Q+\be-\ga)S_b(Q+\al-\ga)}{S_b(\al)S_b(\be)}(1+\CO(e^{2\pi iby}))\\
F_b(\al,\be;\ga;y)=& e^{-2\pi i \al y}
\frac{S_b(\ga)S_b(\al-\be)}{S_b(\be)S_b(\ga-\al)}(1+\CO(e^{-2\pi iby})) \\
&+e^{-2\pi i \be y}
\frac{S_b(\ga)S_b(\be-\al)}{S_b(\al)S_b(\ga-\be)}(1+\CO(-e^{2\pi iby})).
\end{aligned}\end{equation}
\end{lem}
There is also a kind of deformed Euler-integral 
for the hypergeometric function \cite{NU}:
\begin{equation}\label{Eulerint}
\Psi_b(\al,\be;\ga;y)=\frac{1}{i}\int\limits_{-i\infty}^{i\infty}ds \;\,
e^{2\pi i s \be}\frac{G_b(s+y)G_b(s+\ga-\be)}{G_b(s+y+\al)G_b(s+Q)}.
\end{equation}
For the case of main interest, $\Re(\al)=\Re(\be)=Q/2$, $\Re(\ga)=Q$ and
$\Re(x)=0$ one needs to deform the contour such that it passes the pole
at $s=0$ in the right half plane, the pole at $s=-y$ in the
left half plane respectively. It then defines a function that
is analytic in the right $y$ half plane and develops a pole on the
imaginary axis at $x=0$ (Lemma \ref{doubpole}).
\begin{lem} \label{asymeuler}
$\Psi_b(\al,\be;\ga;y)$ has the following asymptotic behavior for
$|\Im(y)|\ra \infty$:
\begin{equation}
\begin{aligned}
\Psi_b(\al,\be;\ga;y)=& 
\frac{G_b(\ga-\be)G_b(\be)}{G_b(\ga)}(1+\CO(e^{2\pi iby})) \\
&+e^{\pi i (\ga-\be)(\ga-\be-Q)}e^{2\pi i(Q-\ga)y}
\frac{G_b(Q+\al-\ga)}{G_b(2Q-\ga)G_b(\al)}(1+\CO(e^{2\pi iby}))\\
\Psi_b(\al,\be;\ga;y)=& e^{-2\pi i \al y}e^{-\pi i \al(\al-Q)}
\frac{G_b(\be-\al)G_b(\ga-\be)}{G_b(\ga-\al)}(1+\CO(e^{-2\pi iby})) \\
&+e^{-2\pi i \be y}e^{-\pi i \be(\be-Q)}
\frac{G_b(\al-\be)G_b(\be)}{G_b(\al)}(1+\CO(e^{-2\pi iby})).
\end{aligned}\end{equation}
\end{lem}
\begin{proof}
In order to study the limit $\Im(y)\ra\infty$
it is convenient to split the integral into two integrals $I_+$ and
$I_-$ over
the intervals $(-y/2,\infty)$ and $(-\infty,-y/2)$ respectively.
In the case of $I_+$ one may use the asymptotics of the $\Theta_b$ functions
containing $y$ for imaginary part of their argument going to $+\infty$, 
eqn. \rf{Thetaas}, to get 
\begin{equation}
\lim_{\Im(y)\ra\infty}I_+=
\lim_{\Im(y)\ra\infty}\frac{1}{i}\int\limits_{-\frac{y}{2}}^{i\infty}
ds \;\,
e^{2\pi i s \be}\frac{G_b(s+\ga-\be)}{G_b(s+Q)}
=\frac{G_b(\be)G_b(\ga-\be)}{G_b(\ga)},
\end{equation}
where \rf{bbeta} was used in the second step.

To study the behavior of $I_-$ for $\Im(y)\ra\infty$ it is convenient to
change the integration variable in the second integral to
$t=s+y$. One gets 
\begin{equation}
I_-=\frac{1}{i}\int\limits_{-i\infty}^{\frac{y}{2}}dt \;\,
e^{2\pi i (t-y) \be}\frac{G_b(t)G_b(t-y+\ga-\be)}{G_b(t+\al)G_b(t-y+Q)}.
\end{equation}
In this expression one may now use
the asymptotics of the $\Theta_b$ functions
containing $y$ for imaginary part of their argument going to $-\infty$,
eqn. \rf{Thetaas}, which yields as previously 
\begin{equation}
\lim_{\Im(y)\ra\infty} e^{-2\pi i y(Q-\ga)}I_-\;=\;
e^{\pi i (\ga-\be)(\ga-\be-Q)}e^{2\pi i(Q-\ga)y}
\frac{G_b(Q+\al-\ga)}{G_b(2Q-\ga)G_b(\al)}.
\end{equation}
The behavior for $\Im(y)\ra-\infty$ is studied similarly.
\end{proof}
\begin{lem}
$\Psi_b(\al,\be;\ga;y)$ is a solution of the finite difference equation
$\CL_b\Psi_b=0$, where
\begin{equation}
\CL_b\equiv
e^{-2\pi i by}(1-T_y^{b})(1-e^{2\pi i b(\ga-Q)}T_y^{b})-
(1-e^{2\pi i b \al}T_y^{b})(1-e^{2\pi i b\be}T_y^{b}).
\end{equation}
\end{lem}
\begin{proof}
Abbreviate the integrand in \rf{Eulerint} by $I$. A direct calculation
shows that it satisfies the equation
\begin{equation}
\CL_b I\;=\; -(1-e^{2\pi i b\al})(1-T_s^{b})
e^{2\pi i s \be}\frac{G_b(s+x)G_b(s+\ga-\be)}{G_b(s+x+\al+b)G_b(s+b^{-1})}.
\end{equation}
The Lemma follows from Cauchy's theorem. 
\end{proof}
The finite difference equation
allows us to define the meromorphic continuation of $\Psi_b$ into 
the right $y$ half plane.
The precise relation between $\Psi_b$ and $F_b$ is
\begin{equation}\label{Psi-F-rel}
\Psi_b(\al,\be;\ga;y)=\frac{G_b(\be)G_b(\ga-\be)}{G_b(\ga)}
F_b(\al,\be;\ga;y'),\qquad y'=y-\fr{1}{2}(\ga-\al-\be+Q).
\end{equation}
This follows as in the proof of Proposition \rf{uniprop} from 
the facts that (i) the finite difference equations satisfied by
left and right hand sides of \rf{Psi-F-rel} are equivalent, and 
(ii) analytic and asymptotic properties of the functions 
of $y$ appearing on both sides of \rf{Psi-F-rel} coincinde.
\newpage
\section{Appendix C}

\theoremstyle{definition}
\newtheorem{Lem}{Lemma}

This appendix collects some results on the analytic and asymptotic properties
of Clebsch-Gordan coefficients, the kernels $\Phi^{\flat}$, $\flat=s,t$
and the Racah-Wigner coefficients.

\subsection{Clebsch-Gordan coefficients}

\begin{Lem}\label{anasC}
The analytic and asymptotic
properties of the Clebsch-Gordan coefficients 
$\CGC{\al_3}{x_3}{\al_2}{x_2}{\al_1}{x_1}$
may be summarized as follows:
\begin{enumerate}
\item $\CGC{Q-\al_3}{x_3}{\al_2}{x_2}{\al_1}{x_1}$ decays 
exponentially as $e^{-2\pi \al_i|x_i|}$ 
if any one of $|x_i|\ra\infty$, $i=1,2,3$.
\item the Clebsch-Gordan coefficients are meromorphic 
w.r.t. each variable $x_i$, $i=1,2,3$ with poles w.r.t. $x_1$ at
\[
\begin{aligned}
\text{Upper half plane: }\qquad & x_1= 
x_2-\fr{i}{2}(\al_1+\al_2-2\al_3)+i(\ep+nb+mb^{-1})\\
 & x_1= x_3-\fr{i}{2}(\al_3+\al_1-Q)+i(\ep+nb+mb^{-1})\\
\text{Lower half plane: }\qquad & x_1= 
x_2-\fr{i}{2}(Q-\al_1-\al_2)-i(Q+nb+mb^{-1})\\
 & x_1= x_3-\fr{i}{2}(2\al_2-\al_3-\al_1)-i(Q+nb+mb^{-1}),
\end{aligned}
\]
where $n,m\in\BZ^{\geq 0}$, and w.r.t. $x_2$ at
\[
\begin{aligned}
\text{Upper half plane: }\qquad & x_2= x_1+\fr{i}{2}(Q-\al_1-\al_2)+i(Q+nb+mb^{-1})\\
 & x_2= 
x_3+\fr{i}{2}(2\al_1-\al_3-\al_2)+i(Q+nb+mb^{-1})\\
\text{Lower half plane: }\qquad & 
x_2= x_1-\fr{i}{2}(2\al_3-\al_1-\al_2)-i(\ep+nb+mb^{-1})\\
& x_2= x_3-\fr{i}{2}(Q-\al_3-\al_2)-i(\ep+nb+mb^{-1}).
\end{aligned}
\]
\end{enumerate}
\end{Lem}
\begin{proof} Direct consequence of analytic and asymptotic properties
of the $S_b$-function given in Appendix B.
\end{proof}
\begin{Lem}\label{anasFTC}
The dependence of $\CGC{\al_3}{\ka_3}{\al_2}{\ka_2}{\al_1}{\ka_1}$
w.r.t. variables $\ka_3$, $\ka_2$, $\ka_1$ is of the following form:
\begin{equation}
\CGC{\al_3}{\ka_3}{\al_2}{\ka_2}{\al_1}{\ka_1}
\;\,=\;\,\de(\ka_3-\ka_2-\ka_1)\;
\CGCZ{\al_3}{\ka_3}{\al_2}{\ka_2}{\al_1}{\ka_1},
\end{equation}
where $\CGCZ{Q-\al_3}{\ka_3}{\al_2}{\ka_2}{\al_1}{\ka_1}$ 
is defined on the hypersurface $\ka_3-\ka_2-\ka_1=0$ only
and is meromorphic w.r.t. $\ka_i$, $i=1,2,3$ with poles only at
\begin{equation}
\ka_i=\pm i(\al_i+nb+mb^{-1}), \qquad i=1,2,3, \qquad n,m\in\BZ^{\geq 0}.
\end{equation}
\end{Lem}
\begin{proof}
One needs to calculate
\begin{equation}
\CGC{\al_3}{\ka_3}{\al_2}{\ka_2}{\al_1}{\ka_1}
\;\,=\;\,
\int\limits_{\BR}dx_2dx_1\; e^{2\pi i k_1x_1}e^{2\pi i k_2x_2}
\CGC{\al_3}{\ka_3}{\al_2}{x_2}{\al_1}{x_1}.
\end{equation}
By inserting \rf{FTk3} and changing variables $(x_1,x_2)\ra (x_+,x_-)$,
$x_{\pm}\equiv x_2\pm x_1$ one finds that the integration over 
$x_+$ produces $\de(\ka_3-\ka_2-\ka_1)$.
$\CGCZ{\al_3}{\ka_3}{\al_2}{\ka_2}{\al_1}{\ka_1}$ is 
therefore given by the
integral
\begin{equation}
\CGCZ{\al_3}{\ka_3}{\al_2}{\ka_2}{\al_1}{\ka_1}\;=\;
\int\limits_{\BR}dx_-\;
e^{\pi i x_-(k_2-k_1)}\;\Phi_{\al_3}
(\al_2,\al_1|\ka_3|x_-).
\end{equation}
It is then useful to employ the Barnes integral representation \rf{Barnesint}
for the b-hypergeometric function that appears in the definition \rf{Phidef}
of the function $\Phi_{\al_3}$. The order of integrals in the resulting
double integral may be exchanged, and the $x_-$ integration carried
out by means of \rf{bbeta}. Up to prefactors that are entire analytic
in $k_i$, $i=1,2,3$ one is left with the following integral:
\begin{equation}
\frac{1}{i}\int\limits_{-i\infty}^{i\infty}ds\;
e^{2\pi i sQ}\frac{G_b(s+A_1)G_b(s+A_2)G_b(s+A_3)}{G_b(s+B_1)
G_b(s+B_2)G_b(s+B_3)},
\end{equation}
where the coefficients are given by
\begin{equation}  
\begin{aligned} A_1=& Q-\al_3+\al_1-\al_2 \\
        A_2=& Q-\al_3-i\ka_3 \\
        A_3=& \al_1+i\ka_1 
\end{aligned}\qquad
\begin{aligned} 
        B_1=& Q+\al_1-\al_2-i\ka_3 \\
        B_2=& 2Q-\al_3-\al_2+i\ka_1 \\ 
        B_3=& Q.
\end{aligned}
\end{equation}
The claim now follows by straightforward application of Lemma \ref{doubpole}.
\end{proof}

\subsection{Kernels $\Phi^{\flat}_{\al_{\flat}}$, $\flat=s,t$}

\begin{lem} \label{anasX}
Analytic and asymptotic properties of 
$\CpBl{\al_s}{\al_3}{\al_2}{\al_4}{\al_1}^{}_{\ep}(x_4;\fx)$
can be summarized
as follows:
\begin{enumerate}
\item 
$\CpBls{\al_s}{\al_3}{\al_2}{\al_4}{\al_1}^{}_{\ep}(x_4;\fx)$ is 
meromorphic w.r.t.
\begin{equation*}\begin{aligned}
x_1 & \quad \text{in}\quad \{x_1\in\BC;\Im(x_1)\in
(-Q,b)\}\\
x_2 & \quad \text{in}\quad \{x_2\in\BC;\Im(x_1)\in
(-b,Q)\}\end{aligned}\qquad
\begin{aligned}
x_3 & \quad \text{in}\quad \{x_3\in\BC;\Im(x_1)\in
(-b,Q)\}\\
x_4 & \quad \text{in}\quad \{x_4\in\BC;\Im(x_1)\in
(-b,b)\}.
\end{aligned}
\end{equation*}
The poles are located at (notation: $x_{ij}\equiv x_i-x_j$)
\begin{equation*}\begin{aligned}
{}& x_{12}+\fr{i}{2}(\al_2+\al_1-2\al_s)-2i\ep=0,\\
{}& x_{12}+\fr{i}{2}(\al_2+\al_1-2(Q-\al_s))-i\ep=0,\\
{}& x_{13}+\fr{i}{2}(\al_3+\al_1-2(Q-\al_4))-2i\ep=0,
\end{aligned}\qquad
\begin{aligned}
{}& x_{14}+\fr{i}{2}(\al_1-\al_4)-2i\ep=0,\\
{}& x_{34}+\fr{i}{2}(\al_4-\al_3)+i\ep=0.
\end{aligned}
\end{equation*}
It decays exponentially for $|x_i|\ra \infty$ as 
$e^{-\pi Q|x_i|}$.
\item $\CpBlt{\al_s}{\al_3}{\al_2}{\al_4}{\al_1}^{}_{\ep}(x_4;\fx)$ is 
analytic w.r.t.
\begin{equation*}\begin{aligned}
x_1 & \quad \text{in}\quad \{x_1\in\BC;\Im(x_1)\in
(-Q,b)\}\\
x_2 & \quad \text{in}\quad \{x_2\in\BC;\Im(x_1)\in
(-Q,b)\}\end{aligned}\qquad
\begin{aligned}
x_3 & \quad \text{in}\quad \{x_3\in\BC;\Im(x_1)\in
(-b,Q)\}\\
x_4 & \quad \text{in}\quad \{x_4\in\BC;\Im(x_1)\in
(-b,b)\}.
\end{aligned}
\end{equation*}
The poles are located at
\begin{equation*}\begin{aligned}
{}& x_{32}-\fr{i}{2}(\al_3+\al_2-2\al_t)+2i\ep=0,\\
{}& x_{32}-\fr{i}{2}(\al_3+\al_2-2(Q-\al_t))+i\ep=0,\\
{}& x_{13}+\fr{i}{2}(\al_3+\al_1-2(Q-\al_4))-2i\ep=0,
\end{aligned}\qquad
\begin{aligned}
{}& x_{14}+\fr{i}{2}(\al_1-\al_4)- i\ep=0,\\
{}& x_{34}+\fr{i}{2}(\al_4-\al_3)+2i\ep=0.
\end{aligned}
\end{equation*}
It decays exponentially for $|x_i|\ra \infty$ as 
$e^{-\pi Q|x_i|}$.
\end{enumerate}
\end{lem}
The residues of these poles that are needed in Section 5 can be 
represented as follows: \begin{equation}
\begin{aligned}
\CR^s_{13}\;\propto & \;\Res_{y_{21}=0}
\CGC{\al_4}{x_4}{\al_3}{x_3}{\al_s}{\ast}\;
\Res_{y_{31}=0}
\CGC{\al_s}{x_s}{\al_2}{x_2}{\al_1}{\ast}_{x_s=
x_3-\frac{i}{2}(\al_s+\al_3-2(Q-\al_4))+i\ep}\\
\CR^s_{14}\;\propto & \;\Res_{y_{31}=0}
\CGC{\al_4}{x_4}{\al_3}{x_3}{\al_s}{\ast}\;\Res_{y_{31}=0}
\CGC{\al_s}{x_s}{\al_2}{x_2}{\al_1}{\ast}_{x_s=
x_4-\frac{i}{2}(\al_s-\al_4)+i\ep}\\
\CR^t_{13}\;\propto & \;\Res_{y_{32}=0}
\CGC{\al_t}{\ast}{\al_3}{x_3}{\al_2}{x_2}\;\Res_{y_{21}=0}
\CGC{\al_4}{x_4'}{\al_t}{x_t}{\al_1}{\ast}_{x_s=
x_3-\frac{i}{2}(\al_3-\al_s)+i\ep}
\\
\CR^t_{14}\;\propto & \;\int_{\BR}dx_t\;\Res_{y_{31}=0}
\CGC{\al_4}{x_4'}{\al_t}{x_t}{\al_1}{\ast}
\CGC{\al_t}{x_t}{\al_3}{x_3}{\al_2}{x_2},
\end{aligned}\end{equation}
where the undetermined prefactor does not depend on any of the variables
and the $\ast$ appearing in the arguments indicates the variable 
of the b-Clebsch-Gordan coefficients that
is to be expressed in terms of the others. 
The necessary residues are
\begin{equation}\begin{aligned}
{} & \Res_{y_{21}=0}\CGC{\al_3}{x_3}{\al_2}{x_2}{\al_1}{\ast}
= \frac{1}{2\pi S_b(\al_3+\al_2+\al_1-Q)}\frac{
S_b\bigl(i(x_3-x_2)-\frac{1}{2}(\al_2-\al_3)\bigr)}{
S_b\bigl(i(x_3-x_2)-\frac{1}{2}(\al_2-\al_3)+\be_{32}\bigr)
}\\
&\qquad\qquad\qquad\qquad\qquad\qquad\qquad\quad
\frac{
S_b\bigl(i(x_2-x_3)+\frac{1}{2}(\al_2+\al_3-2(Q-\al_3))\bigr)}{
S_b\bigl(i(x_2-x_3)+\frac{1}{2}(\al_2+\al_3-2(Q-\al_3))+\be_{31}\bigr)
}\\
{} & \Res_{y_{31}=0}\CGC{\al_3}{\ast}{\al_2}{x_2}{\al_1}{x_1}
= \frac{S_b(\al_3+\al_2-\al_1)}{2\pi}\frac{
S_b\bigl(i(x_1-x_2)-\frac{1}{2}(\al_1+\al_2-2\al_3)\bigr)}{
S_b\bigl(i(x_1-x_2)-\frac{1}{2}(\al_1+\al_2-2\al_3)+\be_{31}\bigr)
}\\
&\qquad\qquad\qquad\qquad\qquad\qquad\qquad\quad
\frac{
S_b\bigl(i(x_1-x_2)-\frac{1}{2}(\al_1+\al_2-2(Q-\al_3))\bigr)}{
S_b\bigl(i(x_1-x_2)-\frac{1}{2}(\al_1+\al_2-2(Q-\al_3))+\be_{32}\bigr)
}\\
{} & \Res_{y_{32}=0}\CGC{\al_3}{x_3}{\al_2}{x_2}{\al_1}{x_1}
= \frac{S_b(\al_3+\al_1-\al_2)}{2\pi}\frac{
S_b\bigl(i(x_1-x_2)-\frac{1}{2}(\al_1+\al_2-2\al_3)\bigr)}{
S_b\bigl(i(x_1-x_2)-\frac{1}{2}(\al_1+\al_2-2\al_3)+\be_{31}\bigr)
}\\
&\qquad\qquad\qquad\qquad\qquad\qquad\qquad\quad
\frac{
S_b\bigl(i(x_1-x_2)-\frac{1}{2}(\al_1+\al_2-2(Q-\al_3))\bigr)}{
S_b\bigl(i(x_1-x_2)-\frac{1}{2}(\al_1+\al_2-2(Q-\al_3))+\be_{21}\bigr)
}\\
& 
\Res_{y_{32}=0}\Res_{y_{21}=0}\CGC{\al_3}{\ast}{\al_2}{\ast}{\al_1}{\ast}
= \Res_{y_{31}=0}\Res_{y_{21}=0}\CGC{\al_3}{\ast}{\al_2}{\ast}{\al_1}{\ast}
=\frac{S_b(2\al_3-Q)}{(2\pi)^2 S_b(\al_1+\al_2+\al_3-Q)}.
\end{aligned}\end{equation}
\begin{lem} \label{anas}
Analytic and asymptotic properties of 
$\CpBl{\al_s}{\al_3}{\al_2}{\al_4}{\al_1}^{}_{\ep}(k_4;\fx)$, $\flat=s,t$
can be summarized
as follows:
\begin{enumerate}
\item 
$\CpBls{\al_s}{\al_3}{\al_2}{\al_4}{\al_1}^{}_{\ep}(k_4;\fx)$ is 
meromorphic w.r.t.
\begin{equation*}\begin{aligned}
x_1 & \;\; \text{in}\;\; \{x_1\in\BC;\Im(x_1)\in
(-Q,b)\},\\
x_2 & \;\; \text{in}\;\; \{x_2\in\BC;\Im(x_1)\in
(-b,Q)\},\end{aligned}\qquad
\begin{aligned}
x_3 &\;\; \text{in}\;\; \{x_3\in\BC;\Im(x_1)\in
(-b,Q)\},\\
k_4 & \;\; \text{in}\;\; \{k_4\in\BC;\Im(x_1)\in
(-\fr{Q}{2},\fr{Q}{2})\}.
\end{aligned}
\end{equation*}
\item $\CpBlt{\al_s}{\al_3}{\al_2}{\al_4}{\al_1}^{}_{\ep}(k_4;\fx)$ is
meromorphic w.r.t.
\begin{equation*}\begin{aligned}
x_1 & \quad \text{in}\quad \{x_1\in\BC;\Im(x_1)\in
(-Q,b)\}\\
x_2 & \quad \text{in}\quad \{x_2\in\BC;\Im(x_1)\in
(-Q,b)\}\end{aligned}\qquad
\begin{aligned}
x_3 & \quad \text{in}\quad \{x_3\in\BC;\Im(x_1)\in
(-b,Q)\}\\
k_4 & \quad \text{in}\quad \{k_4\in\BC;\Im(x_1)\in
(-\fr{Q}{2},\fr{Q}{2})\}.
\end{aligned}
\end{equation*}\end{enumerate}
The poles in their dependence on $x_1,x_2,x_3$ are those poles of 
$\CpBl{\al_s}{\al_3}{\al_2}{\al_4}{\al_1}^{}_{\ep}(x_4;\fx)$, $\flat=s,t$,
which are at positions independent of $x_4$.
Both behave asymptotically
\begin{equation*}\begin{aligned}
\text{for}\;\;|x_1|\ra \infty\;\;& \text{as}\;\; e^{-2\pi i k_4x_1},\\
\text{for}\;\;|x_2|\ra \infty\;\;& \text{as}\;\; e^{-2\pi \al_2|x_2|},
\end{aligned}\qquad
\begin{aligned}
\text{for}\;\;|x_3|\ra \infty\;\; &\text{as}\;\; e^{-2\pi i k_4x_3},\\
\text{for}\;\;|k_4|\ra \infty\;\; & \text{as}\;\; e^{-2\pi \ep k_4}.
\end{aligned}\end{equation*} 
\end{lem}

\subsection{Racah-Wigner coefficients}

\begin{lem}
$\SJS{\al_1}{\al_2}{\al_3}{\al_4}{\al_s}{\al_t}$ is meromorphic w.r.t. 
all six variables and has poles at $\be=-nb-mb^{-1}$ where 
$n,m\in \BZ^{\geq 0}$ and 
$\be$ may be any of the following:
\[
\begin{aligned}
{} &
\begin{aligned}
{}& \al_2+\al_1-\al_s\\
{}& \al_s+\al_1-\al_2
\end{aligned}\quad 
 & 
\begin{aligned}
{}& Q-\al_s-\al_2+\al_1\\
{}& 2Q-\al_1-\al_2-\al_s
\end{aligned} \qquad\qquad
 & 
\begin{aligned}
{}& Q-\al_s-\al_4+\al_3\\
{}& Q-\al_s-\al_3+\al_4
\end{aligned}\quad 
 & 
\begin{aligned}
{}& 2Q-\al_3-\al_4-\al_s\\
{}& Q-\al_3-\al_4+\al_s
\end{aligned}
\\[1ex]
& 
\begin{aligned}
{}& \al_3+\al_2+\al_t-Q\\
{}& \al_3+\al_2-\al_t
\end{aligned}\quad 
 & 
\begin{aligned}
{}& Q-\al_3-\al_t-\al_2\\
{}& Q-\al_2-\al_t-\al_3
\end{aligned} 
\qquad\qquad
 & 
\begin{aligned}
{}& \al_1+\al_4+\al_t-Q\\
{}& \al_1+\al_4-\al_t
\end{aligned}\quad 
 & 
\begin{aligned}
{}& \al_t+\al_4-\al_1\\
{}& Q-\al_1+\al_4-\al_t
\end{aligned}\end{aligned}
\]
\end{lem}
\newcommand{\CMP}[3]{{\it Comm. Math. Phys. }{\bf #1} (#2) #3}
\newcommand{\LMP}[3]{{\it Lett. Math. Phys. }{\bf #1} (#2) #3}
\newcommand{\IMP}[3]{{\it Int. J. Mod. Phys. }{\bf A#1} (#2) #3}
\newcommand{\NP}[3]{{\it Nucl. Phys. }{\bf B#1} (#2) #3}
\newcommand{\PL}[3]{{\it Phys. Lett. }{\bf B#1} (#2) #3}
\newcommand{\MPL}[3]{{\it Mod. Phys. Lett. }{\bf A#1} (#2) #3}
\newcommand{\PRL}[3]{{\it Phys. Rev. Lett. }{\bf #1} (#2) #3}
\newcommand{\AP}[3]{{\it Ann. Phys. (N.Y.) }{\bf #1} (#2) #3}
\newcommand{\LMJ}[3]{{\it Leningrad Math. J. }{\bf #1} (#2) #3}
\newcommand{\FAA}[3]{{\it Funct. Anal. Appl. }{\bf #1} (#2) #3}
\newcommand{\PTPS}[3]{{\it Progr. Theor. Phys. Suppl. }{\bf #1} (#2) #3}
\newcommand{\LMN}[3]{{\it Lecture Notes in Mathematics }{\bf #1} (#2) #2}

\end{document}